\documentclass[11pt]{article}
\usepackage{latexsym,amssymb}
\usepackage{amsmath}
\usepackage{color}

\newtheorem{Theorem}{\bf Theorem}[section]
\newtheorem{Lemma}{\bf Lemma}[section]
\newtheorem{Proposition}{\bf Proposition}[section]
\newtheorem{Corollary}{\bf Corollary}[section]
\newtheorem{Remark}{\bf Remark}[section]
\newtheorem{Example}{\bf Example}[section]
\newtheorem{Definition}{\bf Definition}[section]

\newenvironment{theorem}{\begin{Theorem}$\!\!\!$}{\end{Theorem}}
\newenvironment{lemma}{\begin{Lemma}$\!\!\!$}{\end{Lemma}}
\newenvironment{proposition}{\begin{Proposition}$\!\!\!$}{\end{Proposition}}
\newenvironment{corollary}{\begin{Corollary}$\!\!\!$}{\end{Corollary}}
\newenvironment{remark}{\begin{Remark}$\!\!\!$}{\end{Remark}}

\def\Xint#1{\mathchoice
{\XXint\displaystyle\textstyle{#1}}%
{\XXint\textstyle\scriptstyle{#1}}%
{\XXint\scriptstyle\scriptscriptstyle{#1}}%
{\XXint\scriptscriptstyle\scriptscriptstyle{#1}}%
\!\int}
\def\XXint#1#2#3{{\setbox0=\hbox{$#1{#2#3}{\int}$}
\vcenter{\hbox{$#2#3$}}\kern-.5\wd0}}

\def\dashint{\Xint-}

\numberwithin{equation}{section}

\begin{document}

\title{Existence of solutions for a higher-order semilinear\\ 
parabolic equation with singular initial data}
\author{
\qquad\\
Kazuhiro Ishige, Tatsuki Kawakami and Shinya Okabe
}
\date{}
\maketitle
\begin{abstract}
We establish the existence of solutions 
of the Cauchy problem for a higher-order semilinear parabolic equation 
by introducing a new majorizing kernel. 
We also study necessary conditions on the initial data for the existence of local-in-time solutions 
and identify the strongest singularity of the initial data for the solvability of the Cauchy problem. 
\end{abstract}
\vspace{25pt}
\noindent Addresses:

\smallskip
\noindent K.~I.:  Graduate School of Mathematical Sciences, The University of Tokyo,\\
\qquad\,\,\, 3-8-1 Komaba, Meguro-ku, Tokyo 153-8914, Japan\\
\noindent 
E-mail: {\tt ishige@ms.u-tokyo.ac.jp}\\

\smallskip
\noindent 
T. K.: Department of Applied Mathematics and Informatics, 
Ryukoku University,\\
\qquad\,\,\,\,\, Seta, Otsu 520-2194, Japan\\
\noindent 
E-mail: {\tt kawakami@math.ryukoku.ac.jp}\\

\smallskip
\noindent S. O.:  Mathematical Institute, Tohoku University,
Aoba, Sendai 980-8578, Japan\\
\noindent 
E-mail: {\tt shinya.okabe@tohoku.ac.jp}\\
\vspace{20pt}
%\newline
%\noindent
%{\it 2010 AMS Subject Classifications}: Primary 35B09, 35J61; Secondary 35B32
%\vspace{3pt}
%\newline

\newpage
%%%%%%%%%%%%%%%%%%%%%%%%%%%%%%%%%%%%%
%%%%%%%%%%%%%%%%%%%%%%%%%%%%%%%%%%%%%
\section{Introduction}
%%%%%%%%%%%%%%%%%%%%%%%%%%%%%%%%%%%%%
%%%%%%%%%%%%%%%%%%%%%%%%%%%%%%%%%%%%%
Consider the Cauchy problem 
for a higher-order nonlinear parabolic equation
\begin{equation}
\label{eq:1.1}
\left\{
\begin{array}{ll}
\partial_t u+(-\Delta)^m u=|u|^p,\quad & x\in{\bf R}^N,\,\,t>0,\vspace{3pt}\\
u(x,0)=\mu(x)\ge 0, & x\in{\bf R}^N,
\end{array}
\right.
\end{equation}
where $m=2,3,\dots$, $p>1$ and $\mu$ is a nonnegative measurable function in ${\bf R}^N$ 
or a nonnegative Radon measure in ${\bf R}^N$. 
Problem~\eqref{eq:1.1} is one of the simplest evolution problems 
for higher-order nonlinear parabolic equations. 
In this paper we establish the existence of solutions of problem~\eqref{eq:1.1} 
by introducing a new majorizing kernel to the operator 
\begin{equation}
\label{eq:1.2}
\partial_t u+(-\Delta)^m u=0\quad\mbox{in}\quad{\bf R}^N\times(0,\infty). 
\end{equation}
We also study necessary conditions on the initial data 
for the existence of local-in-time solutions of \eqref{eq:1.1} 
and we identify the strongest singularity of the initial data for the solvability of problem~\eqref{eq:1.1}. 

Before considering problem~\eqref{eq:1.1}, we recall some results on the Cauchy problem 
for a semilinear parabolic equation
\begin{equation}
\label{eq:1.3}
\left\{
\begin{array}{ll}
\partial_t u-\Delta u=u^p,\quad & x\in{\bf R}^N,\,\,t>0,\vspace{3pt}\\
u(x,0)=\mu(x)\ge 0, & x\in{\bf R}^N.
\end{array}
\right.
\end{equation}
In 1985 Baras and Pierre \cite{BP} studied necessary conditions 
for the existence of local-in-time solutions of \eqref{eq:1.3} 
and proved the following (see also \cite{HI01} and \cite{Takahashi}).
\begin{theorem}
\label{Theorem:1.1}
Let $u$ be a nonnegative solution of \eqref{eq:1.3} in ${\bf R}^N\times[0,T)$ 
for some $T>0$, where $\mu$ is a nonnegative Radon measure in ${\bf R}^N$. 
Then there exists $c=c(N,p)>0$ such that 
\begin{equation}
\label{eq:1.4}
\sup_{x\in{\bf R}^N}\mu(B(x,\sigma))\le c\sigma^{N-\frac{2}{p-1}}
\quad\mbox{for}\quad 0<\sigma\le T^{\frac{1}{2}}.
\end{equation}
In particular, in the case of $p=p_1:=1+2/N$, 
there exists $c'=c'(N)>0$ such that
$$
\sup_{x\in{\bf R}^N}\mu(B(x,\sigma))\le 
c\biggr[\log\biggr(e+\frac{T^{\frac{1}{2}}}{\sigma}\biggr)\biggr]^{-\frac{N}{2}}
\quad\mbox{for}\quad 0<\sigma\le T^{\frac{1}{2}}.
$$
\end{theorem}
We remark that, if $1<p<p_1$, then \eqref{eq:1.4} is equivalent to 
\begin{equation}
\label{eq:1.5}
\sup_{x\in{\bf R}^N}\mu(B(x,T^{1/2}))\le cT^{\frac{N}{2}-\frac{1}{p-1}}.
\end{equation}
By Theorem~\ref{Theorem:1.1} we have:
\begin{itemize}
  \item[(a)] 
  There exists $c_1=c_1(N,p)>0$ such that, 
  if $\mu$ is a nonnegative measurable function in ${\bf R}^N$ satisfying 
  $$
  \begin{array}{ll}
  \mu(x)\ge c_1|x|^{-N}\displaystyle{\biggr[\log\left(e+\frac{1}{|x|}\right)\biggr]^{-\frac{N}{2}-1}}
  \quad & \mbox{if}\quad \displaystyle{p=p_1},\vspace{7pt}\\
  \mu(x)\ge c_1|x|^{-\frac{2}{p-1}}\quad & \mbox{if}\quad \displaystyle{p>p_1},
  \end{array}
  $$
  in a neighborhood of the origin, then problem~\eqref{eq:1.3} possesses no local-in-time solutions. 
\end{itemize}
Sufficient conditions for the existence of solutions of problem~\eqref{eq:1.3} 
have been studied in many papers since the pioneering work due to \cite{W1}.  
See e.g. \cite{AD, BaK, BreC, FI, HI01, IKK01, KY, QS, RS, Takahashi, W2} and references therein. 
Among others, 
by \cite{HI01} and \cite{RS} we have: 
\begin{itemize}
\item[(b)] Let $1<p<p_1$. Then there exists $c_2=c_2(N,p)>0$ such that, if
 $$
 \sup_{x\in{\bf R}^N}\mu(B(x,T^{\frac{1}{2}}))\le c_2T^{\frac{N}{2}-\frac{1}{p-1}}
 $$
 for some $T>0$, then problem~\eqref{eq:1.3} possesses a solution in ${\bf R}^N\times[0,T)$.
\item[(c)] Let $p>p_1$. Then there exists $c_3=c_3(N,p)>0$ such that, if 
$$
  \begin{array}{ll}
  0\le\mu(x)\le c_3|x|^{-N}\displaystyle{\biggr[\log\left(e+\frac{1}{|x|}\right)\biggr]^{-\frac{N}{2}-1}}+c_3
  \quad & \mbox{if}\quad \displaystyle{p=p_1},\vspace{7pt}\\
  0\le\mu(x)\le c_3|x|^{-\frac{2}{p-1}}+c_3\quad & \mbox{if}\quad \displaystyle{p>p_1},
  \end{array}
  $$
  then problem~\eqref{eq:1.3} possesses a local-in-time solution.
\end{itemize}
By assertions~(a) and (c) we can identify the strongest singularity of the initial data 
for the existence of solutions of \eqref{eq:1.3} with $p\ge p_1$. 
Assertions~(b) and (c) are proved by the construction of suitable supersolutions 
of \eqref{eq:1.3} and 
the order-preserving property and the semigroup property of the heat operator
are crucial in the proofs. 
\vspace{3pt}

The operator $\partial_t+(-\Delta)^m$ is not order-preserving 
and the study of the solvability of problem~\eqref{eq:1.1} is more delicate than that of problem~\eqref{eq:1.3}. 
Indeed, the fundamental solution $G_m=G_m(x,t)$ of \eqref{eq:1.2} changes its sign for $t>0$. 
In the study of higher-order parabolic equations 
it is crucial to find a suitable majorizing kernel associated with $\partial_t+(-\Delta)^m$. 
Galaktionov and Pohozaev~\cite{GP} found a majorizing kernel of the form 
\begin{equation}
\label{eq:1.6}
\overline{G}_m(x,t):=Dt^{-\frac{N}{2m}}\exp\left(-d\eta^{\frac{2m}{2m-1}}\right)
\quad\mbox{with}\quad \eta=\eta(x,t)=t^{-\frac{1}{2m}}|x|,
\end{equation}
where $D$ and $d$ are positive constants (see Section~2.1), and 
proved the existence of global-in-time solutions of \eqref{eq:1.1} 
for any sufficiently small initial data in $L^1\cap L^\infty$ in the case of $p>p_m:=1+2m/N$. 
They also proved nonexistence of global-in-time solutions of \eqref{eq:1.1} provided that 
$1<p\le p_m$ and $\mu(x)\ge 0$ $(\not\equiv 0)$ in ${\bf R}^N$. 
Subsequently, the existence and the asymptotic behavior of global-in-time solutions with 
bounded initial data have been studied in several papers
under suitable assumptions on the decay of the initial data at the space infinity. 
See e.g. \cite{GP, IKK01, IKK02}. 
(See also \cite{FGG, GG}.) 
On the other hand, 
it does not seem enough to study sufficient conditions 
for the existence of local-in-time solutions 
of problem~\eqref{eq:1.1} with singular initial data, 
although the results in \cite{C} are available. 
As far as we know, 
there are no results related to the identification of the strongest singularity of the initial data 
for the existence of solutions of \eqref{eq:1.1}.  
One of the difficulties is that 
the integral operator associated with $\overline{G}_m$ does not have the semigroup property. 
Indeed, we can not apply 
the arguments in \cite{HI01, RS, W1} with  the majorizing kernel $\overline{G}_m$
to problem~\eqref{eq:1.1}. 
\vspace{3pt}

In this paper, by use of the fundamental solution of 
\begin{equation}
\label{eq:1.7}
\partial_t u+(-\Delta)^{\frac{\theta}{2}}u=0\quad\mbox{in}\quad{\bf R}^N\times(0,\infty),
\end{equation}
where $0<\theta<2$, 
we introduce a new majorizing kernel $K=K(x,t)$ satisfying
\begin{equation}
\label{eq:1.8}
\begin{split}
 & |G_m(x,t)|\le C_1K(x,t),\\
 & \int_{{\bf R}^N}K(x-y,t-s)K(y,s)\,dy\le C_2K(x,t),
\end{split}
\end{equation}
for $x\in{\bf R}^N$ and $0<s<t$. 
Here $C_1$ and $C_2$ are positive constants. 
Applying the arguments in \cite{HI01,TW} with an integral operator associated with $K$, 
we establish the existence of solutions of problem~\eqref{eq:1.1}. 
Furthermore, we modify the arguments in \cite{CM, IkedaSob} to 
study necessary conditions on the initial data for the existence of local-in-time solutions of \eqref{eq:1.1}. 
Then we can identify the strongest singularity of the initial data 
for the existence of local-in-time solutions of \eqref{eq:1.1}.
\vspace{5pt}

Now we are ready to state our main results of this paper. 
The first theorem concerns necessary conditions for the solvability of problem~\eqref{eq:1.1} 
and it corresponds to Theorem~\ref{Theorem:1.1}. 
\begin{theorem}
\label{Theorem:1.2}
Let $N\ge 1$, $m=2,3,\dots$ and $p>1$. 
Let $u$ be a solution of \eqref{eq:1.1} in ${\bf R}^N\times[0,T)$ 
for some $T>0$, where $\mu$ is a nonnegative Radon measure in ${\bf R}^N$. 
Then there exists $\gamma=\gamma(N,m,p)>0$ such that 
\begin{equation}
\label{eq:1.9}
\sup_{x\in{\bf R}^N}\mu(B(x,\sigma))\le \gamma\sigma^{N-\frac{2m}{p-1}}
\quad\mbox{for}\quad 0<\sigma\le T^{\frac{1}{2m}}.
\end{equation}
In particular, if $p=p_m:=1+2m/N$, 
then there exists $\gamma'=\gamma'(N,m)$ such that 
\begin{equation}
\label{eq:1.10}
\sup_{x\in{\bf R}^N}\mu(B(x,\sigma))\le 
\gamma'\biggr[\log\biggr(e+\frac{T^{\frac{1}{2m}}}{\sigma}\biggr)\biggr]^{-\frac{N}{2m}}
\quad\mbox{for}\quad 0<\sigma\le T^{\frac{1}{2m}}.
\end{equation}
\end{theorem} 
Similarly to \eqref{eq:1.5}, 
if $1<p<p_m$, then \eqref{eq:1.9} is equivalent to 
$$
\sup_{x\in{\bf R}^N}\mu(B(x,T^{\frac{1}{2m}}))\le\gamma T^{\frac{N}{2m}-\frac{1}{p-1}}.
$$
As a corollary of Theorem~\ref{Theorem:1.2}, 
we have:
\begin{corollary}
\label{Corollary:1.1}
Let $N\ge 1$, $m=2,3,\dots$ and $p\ge p_m$. 
Then there exists $\gamma_1=\gamma_1(N,m,p)>0$ such that, 
if $\mu$ is a nonnegative measurable function in ${\bf R}^N$ satisfying 
$$
\begin{array}{ll}
\mu(x)\ge \gamma_1|x|^{-N}\displaystyle{\biggr[\log\left(e+\frac{1}{|x|}\right)\biggr]^{-\frac{N}{2m}-1}}
\quad & \mbox{if}\quad \displaystyle{p=p_m},\vspace{7pt}\\
\mu(x)\ge \gamma_1|x|^{-\frac{2m}{p-1}}\quad & \mbox{if}\quad \displaystyle{p>p_m},
\end{array}
$$
in a neighborhood of the origin, then problem~\eqref{eq:1.1} possesses no local-in-time solutions. 
\end{corollary}
Corollary~\ref{Corollary:1.1} corresponds to assertion~(a). 
Next we state results on sufficient conditions 
for the existence of solutions of problem~\eqref{eq:1.1}. 
\begin{theorem}
\label{Theorem:1.3}
Let $N\ge 1$, $m=2,3,\dots$ and $1<p<p_m$. 
Let $\mu$ be a nonnegative Radon measure in ${\bf R}^N$. 
Then there exists $\gamma_2=\gamma_2(N,m,p)>0$ such that, if
\begin{equation}
\label{eq:1.11}
\sup_{x\in{\bf R}^N}\mu(B(x,T^{\frac{1}{2m}}))\le\gamma_2T^{\frac{N}{2m}-\frac{1}{p-1}}
\end{equation}
for some $T>0$, then problem~\eqref{eq:1.1} possesses a solution in ${\bf R}^N\times[0,T)$. 
\end{theorem}
\begin{theorem}
\label{Theorem:1.4}
Let $N\ge 1$, $m=2,3,\dots$ and $p\ge p_m$.
Then there exists $\gamma_3=\gamma_3(N,m,p)>0$ such that, if
$$
\begin{array}{ll}
0\le \mu(x)\le \gamma_3x|^{-N}\displaystyle{\biggr[\log\left(e+\frac{1}{|x|}\right)\biggr]^{-\frac{N}{2m}-1}}+\gamma_3
\quad & \mbox{if}\quad \displaystyle{p=p_m},\vspace{7pt}\\
0\le \mu(x)\le \gamma_3|x|^{-\frac{2m}{p-1}}+\gamma_3\quad & \mbox{if}\quad \displaystyle{p>p_m},
\end{array}
$$
then problem~\eqref{eq:1.1} possesses a local-in-time solution.
\end{theorem}
Theorems~\ref{Theorem:1.3} and \ref{Theorem:1.4} correspond to assertions~(b) and (c), respectively. 
Theorem~\ref{Theorem:1.4} is a direct consequence of Theorems~\ref{Theorem:5.2} and \ref{Theorem:5.3}. 
(See also Remarks~\ref{Remark:5.1} and \ref{Remark:5.2}.)
Furthermore, as a corollary of Theorems~\ref{Theorem:1.2} and \ref{Theorem:1.3}, 
we have:
\begin{corollary}
\label{Corollary:1.2}
Let $\delta$ be the Dirac delta function in ${\bf R}^N$. 
Then problem~\eqref{eq:1.1} possesses a local-in-time solution 
with $\mu=D\delta$ for some $D>0$ if and only if $1<p<p_m$. 
\end{corollary}

The rest of this paper is organized as follows. 
In Section~2 we collect preliminary results on the operators 
$\partial_t+(-\Delta)^m$ $(m=2,3,\dots)$ and $\partial_t+(-\Delta)^{\theta/2}$ $(0<\theta<2)$ 
and their associated semigroups. 
We also formulate the definition of solutions of problem~\eqref{eq:1.1}. 
Furthermore, we formulate the definition of solutions of 
an integral equation associated with problem~\eqref{eq:1.1} 
and prove some properties of the solutions.  
In Section~3 we modify the arguments in \cite{CM,IkedaSob} to prove Theorem~\ref{Theorem:1.2}. 
In Section~4 we introduce a majorizing kernel $K=K(x,t)$ associated with $\partial_t+(-\Delta)^m$ 
and prove \eqref{eq:1.8}. 
In Section~5 we establish the existence of solutions of problem~\eqref{eq:1.1}. 
%%%%%%%%%%%%%%%%%%%%%%%%%%%%%%%%%%%%%
%%%%%%%%%%%%%%%%%%%%%%%%%%%%%%%%%%%%%
\section{Preliminaries}
%%%%%%%%%%%%%%%%%%%%%%%%%%%%%%%%%%%%%
%%%%%%%%%%%%%%%%%%%%%%%%%%%%%%%%%%%%%
This section is divided into three subsections. 
In Sections~2.1 and 2.2 
we recall some preliminary results 
on the operators $\partial_t+(-\Delta)^m$ $(m=2,3,\dots)$ and $\partial_t+(-\Delta)^{\theta/2}$ $(0<\theta<2)$, 
respectively. In Section~2.3 we formulate the definition of solutions of problem~\eqref{eq:1.1}. 
Furthermore, we introduce an integral equation associated with problem~\eqref{eq:1.1} 
and prove some properties of the solutions. 

We introduce some notation. 
For any $1\le r\le\infty$, we denote by $\|\cdot\|_r$ the usual norm of $L^r:=L^r({\bf R}^N)$. 
For any $x\in{\bf R}^N$ and $R>0$, we set $B(x,R):=\{y\in{\bf R}^N\,:\,|x-y|<R\}$. 
For any multi-index $\alpha=(\alpha_1,\dots,\alpha_N)\in {\bf M}:=({\bf N}\cup\{0\})^N$, 
we write 
$$
|\alpha|:=\sum_{i=1}^N\alpha_i,\qquad
\partial_x^\alpha:=\frac{\partial^{|\alpha|}}{\partial x_1^{\alpha_1}\cdots\partial x_N^{\alpha_N}}.
$$
By the letter $C$
we denote generic positive constants and they may have different values also within the same line. 
%%%%%%%
\subsection{Fundamental solutions to $\partial_t+(-\Delta)^m$ $(m=2,3,\dots)$}
%%%%%%%
Let $G_m=G_m(x,t)$ $(m=2,3,\dots)$ be the fundamental solution of \eqref{eq:1.2}. 
Then $G_m$ is represented by 
\begin{equation}
\label{eq:2.1}
G_m(x,t)=(2\pi)^{-\frac{N}{2}}\int_{{\bf R}^N}e^{ix\cdot \xi}e^{-t|\xi|^{2m}}\,d\xi,
\quad
x\in{\bf R}^N,\,\,t>0. 
\end{equation}
The function $G_m$ changes its sign and the operator $\partial_t+(-\Delta)^m$ is not order-preserving.
Let $\overline{G}_m$ be as in \eqref{eq:1.6}. 
Then,
under a suitable choice of $D$ and $d$, it follows that 
\begin{equation}
\label{eq:2.2}
|G_m(x,t)|\le \overline{G}_m(x,t),\quad x\in{\bf R}^N,\,\,\,t>0. 
\end{equation}
(See \cite{GP}.)
Furthermore, $G_m$ satisfies
\begin{align}
\label{eq:2.3}
 & G_m(x,t)=t^{-\frac{N}{2m}}G_m(t^{-\frac{1}{2m}}x,1),\\
\label{eq:2.4}
 & G_m(0,t)>0,\\
\label{eq:2.5}
 & |\partial_x^\alpha G_m(x,t)|\le C_\alpha t^{-\frac{N+|\alpha|}{2m}}\exp\left(-C_\alpha^{-1}\eta^{\frac{2m}{2m-1}}\right)
\quad\mbox{with}\quad\eta=t^{-\frac{1}{2m}}|x|,
\end{align}
for $x\in{\bf R}^N$, $t>0$ and $\alpha\in{\bf M}$, 
where $C_\alpha$ is a positive constant. 
\eqref{eq:2.3} and \eqref{eq:2.4} immediately follow from \eqref{eq:2.1}. 
For \eqref{eq:2.5}, see e.g. \cite[Section~3]{C} and \cite{FGG}.

We define an integral operator associated with $G_m$. 
For any (signed) Radon measure $\mu$ in ${\bf R}^N$, 
we set 
\begin{equation}
\label{eq:2.6}
[S_m(t)\mu](x):=\int_{{\bf R}^N}G_m(x-y,t)\,d\mu(y),\quad x\in{\bf R}^N,\,\,t>0. 
\end{equation}
Similarly, for any measurable function $\phi$ in ${\bf R}^N$, we set 
\begin{equation}
\label{eq:2.7}
[S_m(t)\phi](x):=\int_{{\bf R}^N}G_m(x-y,t)\phi(y)\,dy,\quad x\in{\bf R}^N,\,\,t>0. 
\end{equation}
Let $j=0,1,2,\dots$. 
It follows from the Young inequality and \eqref{eq:2.5} that 
\begin{equation}
\label{eq:2.8}
\|\partial_x^\alpha S_m(t)\phi\|_q\le C_{m}t^{-\frac{N}{2m}\left(\frac{1}{p}-\frac{1}{q}\right)-\frac{j}{2m}}\|\phi\|_p,\quad t>0,
\end{equation}
for $\phi\in L^p$ and $\alpha\in{\bf M}$ with $|\alpha|=j$, 
where $1\le p\le q\le\infty$ and $C_{m}$ is a positive constant independent of $p$ and $q$. 
(See also \cite[Section~2]{C}.) 
Furthermore, 
\begin{equation}
\label{eq:2.9}
\lim_{t\to+0}\|S_m(t)\phi-\phi\|_\infty=0
\end{equation}
for $\phi\in C_0({\bf R}^N)$. 
The convergence rate depends on the modulus of continuity of $\phi$.
%%%%%%%
\subsection{Fundamental solutions to $\partial_t+(-\Delta)^{\theta/2}$ $(0<\theta<2)$}
%%%%%%%
Let $0<\theta<2$. 
Let $G_\theta=G_\theta(x,t)$ be the fundamental solution of \eqref{eq:1.7}, 
that is, 
$$
G_\theta(x,t)=(2\pi)^{-\frac{N}{2}}\int_{{\bf R}^N}e^{ix\cdot \xi}e^{-t|\xi|^\theta}\,d\xi.
$$
Then $G_\theta=G_\theta(x,t)$ is a positive, smooth and radially symmetric function in ${\bf R}^N\times(0,\infty)$ 
and satisfies the following properties (see \cite{BT,BK}): 
\begin{align}
\label{eq:2.10}
& G_\theta(x,t)=t^{-\frac{N}{\theta}}G_\theta(t^{-\frac{1}{\theta}}x,1),\\
\label{eq:2.11}
 &|(\partial_x^\alpha G_\theta)(x,t)|\le C_\alpha t^{-\frac{N+|\alpha|}{\theta}}
\big(1+t^{-\frac{1}{\theta}}|x|\big)^{-N-\theta-|\alpha|},\\
 \label{eq:2.12}
& G_\theta(x,t)\ge 
C t^{-\frac{N}{\theta}}
\big(1+t^{-\frac{1}{\theta}}|x|\big)^{-N-\theta},
\end{align}
for $x\in{\bf R}^N$, $t>0$ and $\alpha\in {\bf M}$,
where $C_\alpha$ is a positive constant.
Furthermore, it follows that 
\begin{equation}
\label{eq:2.13}
G_\theta(x,t)=\int_{{\bf R}^N}G_\theta(x-y,t-s)\,G_\theta(y,s)\,dy,
\quad x\in{\bf R}^N,\,\,0<s<t.
\end{equation}
Similarly to \eqref{eq:2.6} and \eqref{eq:2.7}, 
we set 
$$
[S_\theta(t)\mu](x):=\int_{{\bf R}^N}G_\theta(x-y,t)\,d\mu(y),
\quad
[S_\theta(t)\phi](x):=\int_{{\bf R}^N}G_\theta(x-y,t)\phi(y)\,dy,
$$
for (signed) Radon measure $\mu$ in ${\bf R}^N$ 
and measurable function $\phi$ in ${\bf R}^N$. 
Then, for any $j=0,1,2,\dots$, 
by the Young inequality and \eqref{eq:2.11}
we find $C_j>0$ such that 
$$
\|\partial_x^\alpha S_\theta(t)\phi\|_q\le C_j
t^{-\frac{N}{\theta}(\frac{1}{p}-\frac{1}{q})-\frac{j}{\theta}}\|\phi\|_p,
\qquad t>0,
$$
for $\phi\in L^q$, $1\le p\le q\le\infty$ and $\alpha\in{\bf M}$ with $|\alpha|=j$. 
See e.g. \cite{IKK01}. 
Furthermore, we recall the following lemma on the decay of $\|S_\theta(t)\mu\|_\infty$. 
See \cite[Lemma~2.1]{HI01}.
\begin{lemma}
\label{Lemma:2.1}
Let $\mu$ be a nonnegative Radon measure in ${\bf R}^N$ and $0<\theta<2$. 
Then there exists $C=C(N,\theta)>0$ such that
$$
 \|S_\theta(t)\mu\|_\infty\le Ct^{-\frac{N}{\theta}}\sup_{x\in{\bf R}^N}
 \mu(B(x,t^{\frac{1}{\theta}})),
 \quad t>0. 
$$
\end{lemma}
%%%%%%%
\subsection{Definition of solutions of \eqref{eq:1.1}}
%%%%%%%
We formulate a definition of solutions of problem~\eqref{eq:1.1}. 
\begin{Definition}
\label{Definition:2.1}
Let $N\ge 1$, $m=2,3,\dots$, $p>1$ and $T>0$. 
Let $u$ be a locally integrable function in ${\bf R}^N\times[0,T)$. 
Then we say that $u$ is a solution of \eqref{eq:1.1} in ${\bf R}^N\times[0,T)$ 
if $u$ satisfies 
$$
-\int_{{\bf R}^N}\varphi(x,0)\,d\mu(x)+
\int_0^T\int_{{\bf R}^N}\left[-u\partial_t\varphi+u(-\Delta)^m\varphi\right]\,dx\,dt
=\int_0^T\int_{{\bf R}^N}|u|^p\varphi\,dx\,dt
$$
for all $\varphi\in C^\infty({\bf R}^N\times[0,T))$ with $\mbox{{\rm supp}}\,\varphi\subset B(0,R)\times[0,T-\epsilon]$ 
for some $R>0$ and $0<\epsilon<T$. 
\end{Definition}
We also formulate a definition of solutions of the integral equation 
\begin{equation}
\tag{I}
u(x,t)=\int_{{\bf R}^N}G_m(x-y,t)\,d\mu(y)+\int_0^t\int_{{\bf R}^N}G_m(x-y,t-s)|u(y,s)|^p\,dy\,ds.
\end{equation}
\begin{Definition}
\label{Definition:2.2}
Let $N\ge 1$, $m=2,3,\dots$, $p>1$ and $\mu$ be a nonnegative Radon measure in ${\bf R}^N$. 
Let $u$ be a continuous function in ${\bf R}^N\times(0,T)$ for some $T>0$ 
and set 
\begin{equation}
\label{eq:2.14}
\begin{split}
 & \overline{u}_1(x,t):=\int_{{\bf R}^N}|G_m(x-y,t)|\,d\mu(y),\\
 & \overline{u}_2(x,t):=\int_0^t\int_{{\bf R}^N}|G_m(x-y,t-s)||u(y,s)|^p\,dy\,ds. 
\end{split}
\end{equation}
We say that $u$ is a solution of integral equation~{\rm (I)} in ${\bf R}^N\times[0,T)$ 
if 
\begin{equation}
\label{eq:2.15}
\sup_{\tau\le t<T}\|\overline{u}_1(t)\|_\infty+\sup_{\tau\le t<T}\|\overline{u}_2(t)\|_\infty<\infty
\quad\mbox{for}\quad \tau\in(0,T)
\end{equation}
and $u$ satisfies integral equation~{\rm (I)} for $(x,t)\in{\bf R}^N\times(0,T)$. 
\end{Definition}
In the rest of this section 
we show that the solution of integral equation~(I) is a solution of \eqref{eq:1.1}. 
\begin{proposition}
\label{Proposition:2.1}
Let $u$ be a solution of integral solution~{\rm (I)} in ${\bf R}^N\times[0,T)$ for some $T>0$. 
\begin{itemize}
  \item[{\rm (a)}]
  For any $\tau\in(0,T)$, 
  $u_\tau$ defined by $u_\tau(x,t):=u(x,t+\tau)$ is a solution of problem~\eqref{eq:1.1} 
  in ${\bf R}^N\times[0,T-\tau)$ with the initial data $u(\tau)$. 
  \item[{\rm (b)}] 
  Let $\alpha\in{\bf M}$ and $i\in\{0,1\}$ be such that $|\alpha|+4i\le 2m$. 
  Then $\partial_t^i\partial_x^\alpha u\in BC({\bf R}^N\times[\tau,T))$ for $\tau\in(0,T)$.
  \item[{\rm (c)}]  
  $u$ satisfies 
  \begin{equation}
  \label{eq:2.16}
  \partial_t u+(-\Delta)^m u=|u|^p,\quad (x,t)\in{\bf R}^N\times(0,T),
  \end{equation}
  in the classical sense. 
\end{itemize}
Furthermore, 
$u$ is a solution of \eqref{eq:1.1} in ${\bf R}^N\times[0,T)$. 
\end{proposition}
{\bf Proof of assertions~(a), (b) and (c).} 
Let $u$ be a solution of integral equation~(I) in ${\bf R}^N\times[0,T)$ for some $T>0$. 
By \eqref{eq:1.6}, \eqref{eq:2.2} and \eqref{eq:2.15} we see that 
\begin{equation*}
\begin{split}
 & \int_{{\bf R}^N}|G_m(x-y,t-\tau)|
\biggr[\int_{{\bf R}^N}|G_m(y-z,\tau)|\,d\mu(z)\biggr]\,dy<\infty,\\
 & \int_{{\bf R}^N}|G_m(x-y,t-\tau)|
\biggr
[\int_0^\tau\int_{{\bf R}^N}|G_m(y-z,\tau-s)|u(z,s)|^p\,dz\,ds\biggr]\,dy<\infty,
\end{split}
\end{equation*}
for $x\in{\bf R}^N$ and $0<\tau<t$. 
It follows from the Fubini theorem that 
\begin{equation*}
\begin{split}
 & \int_{{\bf R}^N}G_m(x-y,t-\tau)u(y,\tau)\,dy\\
 & =\int_{{\bf R}^N}G_m(x-y,t-\tau)\\
 & \qquad\times\biggr
[\int_{{\bf R}^N}G_m(y-z,\tau)\,d\mu(z)+\int_0^\tau\int_{{\bf R}^N}G_m(y-z,\tau-s)|u(z,s)|^p\,dz\,ds\biggr]\,dy\\
 & =\int_{{\bf R}^N}\biggr(\int_{{\bf R}^N}G_m(x-y,t-\tau)G_m(y-z,\tau)\,dy\biggr)\,d\mu(z)\\
 & \qquad+\int_0^\tau\int_{{\bf R}^N}\biggr(\int_{{\bf R}^N}G_m(x-y,t-\tau)G_m(y-z,\tau-s)\,dy\biggr)|u(z,s)|^p\,dz\,ds\\
 & =\int_{{\bf R}^N}G_m(x-z,t)\,d\mu(z)+\int_0^\tau\int_{{\bf R}^N}G_m(x-z,t-s)|u(z,s)|^p\,dz\,ds
\end{split}
\end{equation*}
for $x\in{\bf R}^N$ and $0<\tau<t$. 
This together with Definition~\ref{Definition:2.2} 
implies that 
\begin{align}
\label{eq:2.17}
 & \sup_{\tau\le t<T}\|u(t)\|_\infty<\infty,\quad 0<\tau<T,\\
\notag
 & u(x,t)=\int_{{\bf R}^N}G_m(x-y,t-\tau)u(y,\tau)\,dy\\
\notag
 & \qquad\qquad\quad
+\int_\tau^t\int_{{\bf R}^N}G_m(x-y,t-s)|u(y,s)|^p\,dy\,ds,\quad x\in{\bf R}^N,\,\,0<\tau<T,
\end{align}
and assertion~(a) holds.
%This implies assertion~(a). 
By \eqref{eq:2.17} we apply similar arguments in regularity theorems  
for second order parabolic equations (see e.g. \cite[Chapter 1]{F}) to integral equation~(I)
and 
%we 
obtain assertions~(b) and (c). 
$\Box$
\vspace{3pt}
\newline
It remains to prove that $u$ is a solution of problem~\eqref{eq:1.1}. 
For this aim, we modify the arguments in \cite{HI01} to prepare the following two lemmas. 
\begin{lemma}
\label{Lemma:2.2}
Let $u$ be a solution of integral equation~{\rm (I)} in ${\bf R}^N\times[0,T)$ for some $T>0$. 
Then 
\begin{eqnarray}
\label{eq:2.18}
 & & \lim_{t\to +0}\int_{{\bf R}^N\setminus B(0,R)}\overline{G}_m(\lambda x,t)\,d\mu(x)=0,\\
\label{eq:2.19}
 & & \lim_{t\to +0}\int_0^t\int_{{\bf R}^N\setminus B(0,R)}\overline{G}_m(\lambda x,t-s)|u(x,s)|^p\,dx\,ds=0,\qquad
\end{eqnarray}
for $R>0$ and $\lambda>0$.
\end{lemma}
{\bf Proof.}
By \eqref{eq:2.4} 
we find $R_*>0$ and $c_*>0$ such that 
$$
\inf_{x\in B(0,R_*)}G_m(x,1)\ge c_*>0. 
$$
Then it follows from \eqref{eq:2.3} that
$$
G_m(x-y,t)\ge t^{-\frac{N}{2m}}c_*\quad\mbox{for}\quad x-y\in B(0,R_*t^{\frac{1}{2m}}).
$$
This together with \eqref{eq:2.14} and \eqref{eq:2.15} implies that
\begin{equation*}
\begin{split}
\infty>\|\overline{u}_1(T_\epsilon)\|_\infty & \ge{\overline{u}_1}(x,T_\epsilon)\ge
\int_{B(x,T^{\frac{1}{2m}}R_*)}|G_m(x-y,T_\epsilon)|\,d\mu(y)\\
 & \ge c_*T^{-\frac{N}{2m}}\mu(B(x,R_*T_\epsilon^{\frac{1}{2m}})),\\
\infty>\|\overline{u}_2(T_\epsilon)\|_\infty & \ge
\overline{u}_2(x,T_\epsilon)\ge\int_0^{T_{2\epsilon}}
\int_{B(x,R_*(T_\epsilon-s)^{\frac{1}{2m}})}|G_m(x-y,T_\epsilon-s)|u(y,s)|^p\,dy\,ds\\
 & \ge c_*\int_0^{T_{2\epsilon}}(T_\epsilon-s)^{-\frac{N}{2m}}\int_{B(x,R_*(T_\epsilon-s)^{\frac{1}{2m}})}|u(y,s)|^p\,dy\,ds\\
 & \ge c_*\epsilon^{-\frac{N}{2m}}\int_0^{T_{2\epsilon}}\int_{B(x,\epsilon^{\frac{1}{2m}}R_*)}|u(y,s)|^p\,dy\,ds,
\end{split}
\end{equation*}
for $x\in{\bf R}^N$, where $T_\epsilon:=T-\epsilon$, $T_{2\epsilon}=T-2\epsilon$ and $0<2\epsilon<T$. 
Since $x\in{\bf R}^N$ is arbitrary, 
we deduce that 
\begin{equation}
\label{eq:2.20}
\begin{split}
 & \sup_{x\in{\bf R}^N}\mu(B(x,R))<\infty,\\
 & \sup_{x\in{\bf R}^N}\int_0^{T-\epsilon}\int_{B(x,R)}|u(y,s)|^p\,dy\,ds<\infty,
\end{split}
\end{equation}
for $R>0$ and $0<\epsilon<T/2$. (See \cite[Lemma~2.1]{IS}.)

Let $0<R<\infty$ and set $R':=\min\{R/2,1/2\}$. 
By the Besicovitch covering lemma we can find an integer $n_*$
depending only on $N$ and a set $\{x_{k,i}\}_{k=1,\dots,n_*,\,i\in{\bf N}}\subset{\bf R}^N\setminus B(0,R)$ such that 
\begin{equation}
\label{eq:2.21}
\begin{split}
 & \overline{B(x_{k,i},R')}\cap\overline{B(x_{k,j},R')}\quad\mbox{if}\quad i\not=j,\\
 & {\bf R}^N\setminus B(0,R)\subset\bigcup_{k=1}^{n_*}\bigcup_{i=1}^\infty \overline{B(x_{k,i},R')}
\subset {\bf R}^N\setminus B(0,R/2).
\end{split}
\end{equation}
Then we have 
\begin{equation}
\label{eq:2.22}
\begin{split}
\int_{{\bf R}^N\setminus B(0,R)}\overline{G}_m(\lambda x,t)\,d\mu(x)
 & \le\sum_{k=1}^{n_*}\sum_{i=1}^\infty\int_{\overline{B(x_{k,i},R')}}\overline{G}_m(\lambda x,t)\,d\mu(x)\\
 & \le\sum_{k=1}^{n_*}\sum_{i=1}^\infty\mu(\overline{B(x_{k,i},R')})\sup_{x\in\overline{B(x_{k,i},R')}}\overline{G}_m(\lambda x,t)\\
 & \le\sup_{x\in {\bf R}^N}\mu(B(x,1))\sum_{k=1}^{n_*}\sum_{i=1}^\infty\sup_{x\in\overline{B(x_{k,i},R')}}\overline{G}_m(\lambda x,t).
\end{split}
\end{equation}
Let $\epsilon>0$ be such that $2(1-\epsilon)>1+\epsilon$.
For $k=1,\dots,n_*$ and $i\in \bf N$, 
since $x_{k,i}\not\in B(0,R)$ and $R'\le R/2$, 
we have 
$$
\frac{|x_{k,i}|}{R'}\ge\frac{R}{R'}\ge2>\frac{1+\epsilon}{1-\epsilon},
$$
which implies that $|x_{k,i}|-R'\ge\epsilon(|x_{k,i}|+R')$. 
Then it holds that
$$
|y|\ge |x_{k,i}|-R'\ge\epsilon(|x_{k,i}|+R')\ge\epsilon|z|
$$
for $y$, $z\in\overline{B(x_{k,i},R')}$, $k=1,\dots,n_*$ and $i\in\bf N$. 
Therefore we observe from \eqref{eq:1.6} that 
$$
\sup_{x\in\overline{B(x_{k,i},R')}}\overline{G}_m(\lambda x,t)
\le\inf_{x\in\overline{B(x_{k,i},R')}}\overline{G}_m(\lambda\epsilon x,t)
\le\frac{1}{|B(0,R')|}\int_{\overline{B(x_{k,i},R')}}\overline{G}_m(\lambda\epsilon z,t)\,dz
$$
for $k=1,\dots,n_*$ and $i\in\bf N$, 
where $|B(0,R')|$ is the volume of $B(0,R')$. 
This together with \eqref{eq:2.21} implies that 
\begin{equation}
\label{eq:2.23}
\begin{split}
\sum_{k=1}^{n_*}\sum_{i=1}^\infty\sup_{x\in\overline{B(x_{k,i},R')}}\overline{G}_m(\lambda x,t)
 & \le Cn_*R'^{-N}\int_{{\bf R}^N\setminus B(0,R/2)}\overline{G}_m(\lambda\epsilon z,t)\,dz\\
 & =Cn_*R'^{-N}\int_{{\bf R}^N\setminus t^{-\frac{1}{2m}}B(0,R/2)}\overline{G}_m(\lambda\epsilon z,1)\,dz\to 0
\end{split}
\end{equation}
as $t\to+0$. 
Combining \eqref{eq:2.22} and \eqref{eq:2.23}, 
we obtain \eqref{eq:2.18}. 

Since 
\begin{equation*}
\begin{split}
\overline{G}_m(\lambda x,t-s)
 & \le C(t-s)^{-\frac{N}{2m}}\exp\left(-C^{-1}\eta(\lambda x,t-s)^{\frac{2m}{2m-1}}\right)\\
 & \le \exp\left(-(2C)^{-1}\eta(\lambda x,t-s)^{\frac{2m}{2m-1}}\right)\\
 & \le \exp\left(-(2C)^{-1}\eta(\lambda x,t)^{\frac{2m}{2m-1}}\right)=:\hat{G}_m(\lambda x,t)
\end{split}
\end{equation*}
for $x\in{\bf R}^N\setminus B(0,R)$ and $0<s<t$, 
we have 
\begin{equation}
\label{eq:2.24}
\begin{split}
 & \int_0^t\int_{{\bf R}^N\setminus B(0,R)}\overline{G}_m(\lambda x,t-s)|u(x,s)|^p\,dx\,ds\\
 & \le\sum_{k=1}^{n_*}\sum_{i=1}^\infty\int_0^t\int_{\overline{B(x_{k,i},R')}}\overline{G}_m(\lambda x,t-s)|u(x,s)|^p\,dx\,ds\\
 & \le\sum_{k=1}^{n_*}\sum_{i=1}^\infty\sup_{x\in\overline{B(x_{k,i},R')}}\hat{G}_m(\lambda x,t)
 \int_0^t\int_{\overline{B(x_{k,i},R')}}|u(x,s)|^p\,dx\,ds\\
 & \le\sup_{x\in {\bf R}^N}\int_0^{T/2}\int_{B(x,1)}|u(x,s)|^p\,dx\,ds\,
 \sum_{k=1}^{n_*}\sum_{i=1}^\infty\sup_{x\in\overline{B(x_{k,i},R')}}\hat{G}_m(\lambda x,t)
\end{split}
\end{equation}
for $0<t\le T/2$.
Similarly to \eqref{eq:2.23}, we observe that 
\begin{equation}
\label{eq:2.25}
\sum_{k=1}^{n_*}\sum_{i=1}^\infty\sup_{x\in\overline{B(x_{k,i},R')}}\hat{G}_m(\lambda x,t)
\le Cn_*R'^{-N}\int_{{\bf R}^N\setminus B(0,R/2)}\hat{G}_m(\lambda\epsilon z,t)\,dz
\to 0
\end{equation}
as $t\to+0$. 
Combining \eqref{eq:2.24} and \eqref{eq:2.25}, we see that 
$$
\lim_{t\to +0}\int_0^t\int_{{\bf R}^N\setminus B(0,R)}\overline{G}_m(\lambda x,t-s)|u(x,s)|^p\,dx\,ds=0,
$$
which implies \eqref{eq:2.19}. Thus Lemma~\ref{Lemma:2.2} follows. 
$\Box$
\begin{lemma}
\label{Lemma:2.3}
Let $u$ be a solution of integral equation~{\rm (I)} in ${\bf R}^N\times[0,T)$ for some $T>0$. 
Then 
\begin{eqnarray}
\label{eq:2.26}
 & & 
\lim_{t\to+0}
\int_{{\bf R}^N}\int_{{\bf R}^N}\varphi(x,t)G_m(x-y,t)\,d\mu(y)\,dx
=\int_{{\bf R}^N}\varphi(y,0)\,d\mu(y),\qquad\\
\label{eq:2.27}
 & & \lim_{t\to +0}\int_{{\bf R}^N}\int_0^t\int_{{\bf R}^N}G_m(x-y,t-s)\varphi(x,t)|u(y,s)|^p\,dy\,ds\,dx=0,
\end{eqnarray}
for $\varphi\in C^\infty({\bf R}^N\times[0,T))$ with $\mbox{{\rm supp}}\,\varphi\subset B(0,R)\times[0,T-\epsilon]$ 
for some $R>0$ and $\epsilon\in(0,T)$.
\end{lemma}
{\bf Proof.}
Let $\varphi\in C^\infty({\bf R}^N\times[0,T))$ be such that $\mbox{{\rm supp}}\,\varphi\subset B(0,R)\times[0,T-\epsilon]$ 
for some $R>0$ and $\epsilon\in(0,T)$.
Set 
$$
\Phi(x,t:\tau):=[S_m(t)\varphi(\tau)](x)=\int_{{\bf R}^N}G_m(x-y,t)\varphi(y,\tau)\,dy,
\quad x\in{\bf R}^N,\,\,\,t>0,\,\,\,\tau\in(0,T).
$$
By \eqref{eq:2.8} we have 
\begin{equation}
\label{eq:2.28}
\|\Phi(t:\tau)\|_\infty\le C\|\varphi(\tau)\|_\infty\le C\|\varphi\|_{L^\infty({\bf R}^N\times(0,T))},
\quad t>0,\,\,\,\tau\in(0,T). 
\end{equation}
On the other hand, 
it follows from the Fubini theorem that
\begin{equation}
\label{eq:2.29}
\begin{split}
 & \int_{{\bf R}^N}\int_{{\bf R}^N}\varphi(x,t)G_m(x-y,t)\,d\mu(y)\,dx\\
 & =\int_{{\bf R}^N}\int_{{\bf R}^N}\varphi(x,t)G_m(x-y,t)\,dx\,d\mu(y)
 =\int_{{\bf R}^N}\int_{{\bf R}^N}\varphi(x,t)G_m(y-x,t)\,dx\,d\mu(y)\\
 & =\int_{{\bf R}^N}\Phi(y,t;t)\,d\mu(y)
 =\int_{{\bf R}^N}\varphi(y,0)\,d\mu(y)+\int_{{\bf R}^N}[\Phi(y,t:t)-\varphi(y,0)]\,d\mu(y).
\end{split}
\end{equation}
Since $|x-y|\ge |x|/2$ for $x\in{\bf R}^N\setminus B(0,2R)$ and $y\in B(0,R)$,
by \eqref{eq:2.2}
we can find $\lambda>0$ such that
\begin{equation}
\label{eq:2.30}
\begin{split}
|\Phi(x,t;\tau)| & \le\|\varphi\|_{L^\infty({\bf R}^N\times(0,T))}\int_{B(0,R)}|G_m(x-y,t)|\,dy\\
 & \le C\|\varphi\|_{L^\infty({\bf R}^N\times(0,T))}\overline{G}_m(\lambda x,t)
\end{split}
\end{equation}
for $x\in{\bf R}^N\setminus B(0,2R)$, $t>0$ and $\tau\in(0,T)$. 
Furthermore, 
by the uniform continuity of $\varphi$ in $\overline{B(0,2R)}\times[0,T-\epsilon]$ and \eqref{eq:2.9}
we observe that
\begin{equation}
\label{eq:2.31}
\begin{split}
 & \sup_{x\in B(0,2R)}|\Phi(x,t:t)-\varphi(x,0)|\\
 & \le\sup_{x\in B(0,2R)}|\Phi(x,t:t)-\varphi(x,t)|+\sup_{x\in B(0,2R)}|\varphi(x,t)-\varphi(x,0)|
 \to 0
\end{split}
\end{equation}
as $t\to +0$. 
Therefore, by \eqref{eq:2.18}, \eqref{eq:2.28} and \eqref{eq:2.31} 
we have
\begin{equation}
\label{eq:2.32}
\begin{split}
 & \biggr|\int_{{\bf R}^N}[\Phi(y,t:t)-\varphi(y,0)]\,d\mu(y)\biggr|\\
 & \le\int_{B(0,2R)}|\Phi(y,t:t)-\varphi(y)|\,d\mu(y)
 +\int_{{\bf R}^N\setminus B(0,2R))}|\Phi(y,t:t)|\,d\mu(y)\\
 & \le\sup_{x\in B(0,2R)}|\Phi(x,t:t)-\varphi(x,0)\,|\mu(B(0,2R))\\
 & \qquad\qquad
 +C\|\varphi\|_{L^\infty({\bf R}^N\times(0,T))}\int_{{\bf R}^N\setminus B(0,2R)}\overline{G}_m(\lambda y,t)\,d\mu(y)
 \to 0
\end{split}
\end{equation}
as $t\to +0$. 
Combining \eqref{eq:2.29} and \eqref{eq:2.32}, 
we have \eqref{eq:2.26}. 
Furthermore, 
by \eqref{eq:2.19}, \eqref{eq:2.28} and \eqref{eq:2.30} we have
\begin{equation*}
\begin{split}
 & \left|\int_{{\bf R}^N}\int_0^t\int_{{\bf R}^N}G_m(x-y,t-s)\varphi(x,t)|u(y,s)|^p\,dy\,ds\,dx\right|\\
 & =\left|\int_0^t\int_{{\bf R}^N}\Phi(y,t-s:t)|u(y,s)|^p\,dy\,ds\right|\\
 & \le C\|\varphi\|_{L^\infty({\bf R}^N\times(0,T))}\int_0^t\int_{B(0,2R)}|u(y,s)|^p\,dy\,ds\\
 & \qquad\qquad
 +C\|\varphi\|_{L^\infty({\bf R}^N\times(0,T))}\int_0^t\int_{{\bf R}^N\setminus B(0,2R)}
 \overline{G}_m(\lambda y,t-s)|u(y,s)|^p\,dy\,ds\to 0
\end{split}
\end{equation*}
as $t\to +0$. This implies \eqref{eq:2.27}. 
Thus Lemma~\ref{Lemma:2.3} follows.
$\Box$
\vspace{5pt}

Now we are ready to complete the proof of Proposition~\ref{Proposition:2.1}.
\vspace{5pt}
\newline
{\bf Proof of  Proposition~\ref{Proposition:2.1}.}
Let $u$ be a solution of integral equation~(I) in ${\bf R}^N\times[0,T)$ for some $T>0$. 
It suffices to prove that $u$ is a solution of \eqref{eq:1.1} in ${\bf R}^N\times[0,T)$. 

Let $\varphi\in C^\infty({\bf R}^N\times[0,T))$ be such that 
$\mbox{{\rm supp}}\,\varphi\subset B(0,R)\times[0,T-\epsilon]$ 
for some $R>0$ and $\epsilon\in(0,T)$.
Then it follows from Definition~\ref{Definition:2.2} and Lemma~\ref{Lemma:2.3} that 
\begin{equation}
\label{eq:2.33}
\begin{split}
 & \int_{{\bf R}^N}u(x,t)\varphi(x,t)\,dx\\
 & =\int_{{\bf R}^N}\int_{{\bf R}^N}\varphi(x,t)G_m(x-y,t)\,d\mu(y)\,dx\\
 & +\int_{{\bf R}^N}\int_0^t\int_{{\bf R}^N}\varphi(x,t)G_m(x-y,t-s)|u(y,s)|^p\,dy\,ds\,dx
 \to\int_{{\bf R}^N}\varphi(x,0)\,d\mu(x)
\end{split}
\end{equation}
as $t\to+0$. 
On the other hand, 
by \eqref{eq:2.16} we see that 
$$
-\int_{{\bf R}^N}\varphi(x,\tau)u(x,\tau)\,dx+
\int_\tau^T\int_{{\bf R}^N}\left[-u\partial_t\varphi+u(-\Delta)^m\varphi\right]\,dx\,dt
=\int_\tau^T\int_{{\bf R}^N}|u|^p\varphi\,dx\,dt.
$$
Letting $\tau\to +0$, 
by \eqref{eq:2.20} and \eqref{eq:2.33} we have
$$
-\int_{{\bf R}^N}\varphi(x,0)\,d\mu(x)+
\int_0^T\int_{{\bf R}^N}\left[-u\partial_t\varphi+u(-\Delta)^m\varphi\right]\,dx\,dt
=\int_0^T\int_{{\bf R}^N}|u|^p\varphi\,dx\,dt.
$$
This means that $u$ is a solution of \eqref{eq:1.1} in ${\bf R}^N\times[0,T)$. 
Thus Proposition~\ref{Proposition:2.1} follows.
$\Box$
%%%%%%%%%%%%%%%%%%%%%%%%%%%%%%%%%%%%%
%%%%%%%%%%%%%%%%%%%%%%%%%%%%%%%%%%%%%
\section{Proof of Theorem~\ref{Theorem:1.2}}
%%%%%%%%%%%%%%%%%%%%%%%%%%%%%%%%%%%%%
%%%%%%%%%%%%%%%%%%%%%%%%%%%%%%%%%%%%%
In this section we modify the arguments in \cite{IkedaSob} (see also \cite{CM}) 
to prove Theorem~\ref{Theorem:1.2}.
\vspace{3pt}
\newline
{\bf Proof of Theorem~\ref{Theorem:1.2}.} 
Let $u$ be a solution of problem~\eqref{eq:1.1} in ${\bf R}^N\times[0,T)$ 
for some $T>0$. Set 
\begin{equation}
\label{eq:3.1}
u_T(x,t):=T^{\frac{1}{p-1}}u(T^{\frac{1}{2m}}x,Tt),\qquad
\mu_T(x):=T^{\frac{1}{p-1}}\mu(T^{\frac{1}{2m}}x).
\end{equation}
Then $u_T$ is a solution of problem~\eqref{eq:1.1} in ${\bf R}^N\times[0,1)$ 
with the initial data $\mu_T$. 
Due to similar transformation~\eqref{eq:3.1}, 
it suffices to consider the case of $T=1$ for the proof of Theorem~\ref{Theorem:1.2}. 

Let 
$$
f(s):=e^{-\frac{1}{s}}\quad\mbox{if}\quad s>0,
\qquad 
f(s)=0\quad\mbox{if}\quad s\le 0.
$$
Set 
$$
\eta(s):=\frac{f(2-s)}{f(2-s)+f(s-1)}.
$$
Then $\eta\in C^\infty([0,\infty))$ and 
\begin{equation*}
\begin{split}
 & \eta'(s)=\frac{-f'(2-s)f(s-1)-f(2-s)f'(s-1)}{[f(2-s)+f(s-1)]^2}\le 0,\qquad s\ge 0,\\
 &  \eta(s)=1\quad\mbox{on}\quad[0,1],\qquad \eta(s)=0\quad\mbox{on}\quad[2,\infty).
\end{split}
\end{equation*}
Set 
$$
\eta^*(s)=0\quad\mbox{for}\quad s\in[0,1),
\qquad
\eta^*(s)=\eta(s)\quad\mbox{for}\quad s\ge 1. 
$$
Since $p>1$, 
for $k=1,2,\dots$, it follows that
\begin{equation}
\label{eq:3.2}
|\eta^{(k)}(s)|\le C\eta^*(s)^{\frac{1}{p}}\quad\mbox{for}\quad s\ge 1. 
\end{equation}

Let $u$ be a solution of problem~\eqref{eq:1.1} in ${\bf R}^N\times[0,1)$. 
Let $x_0\in {\bf R}^N$ and $0<r_*<1$ be such that 
$$
\mu\biggr(B\left(x_0,(r_*/3)^{\frac{1}{2m}}\right)\biggr)>0.
$$
For any $R\in(0,1]$, we set
$$
\psi_R(x,t):=\eta\left(3\frac{|x-x_0|^{2m}+t}{R}\right),
\qquad	
\psi_R^*(x,t):=\eta^*\left(3\frac{|x-x_0|^{2m}+t}{R}\right). 
$$
By \eqref{eq:3.2}, for $k=1,2,\dots$, 
we have 
\begin{equation}
\label{eq:3.3}
|\partial_t\psi_R(x,t)|\le CR^{-1}\psi_R^*(x,t)^{\frac{1}{p}},
\qquad
|\nabla_x^k\psi_R(x,t)|\le CR^{-\frac{k}{2m}}\psi_R^*(x,t)^{\frac{1}{p}},
\end{equation}
for $x\in{\bf R}^N$ and $0<t\le 1$. 
It follows from \eqref{eq:3.3} that 
\begin{equation}
\label{eq:3.4}
\begin{split}
 & \int_{{\bf R}^N}\psi_R(x,0)\,d\mu+\int_0^R\int_{{\bf R}^N}|u(x,t)|^p\psi_R(x,t)\,dx\,dt\\
 & =\int_0^R\int_{{\bf R}^N}u(x,t)(-\partial_t+(-\Delta)^m)\psi_R(x,t)\,dx\,dt\\
 & \le CR^{-1}\int_0^R\int_{{\bf R}^N}|u(x,t)|\psi_R^*(x,t)^{\frac{1}{p}}\,dx\,dt\\
 & \le CR^{-1}\left(\int_0^R\int_{{\bf R}^N}\chi_{\{\psi_R^*(x,t)>0\}}\,dx\,dt\right)^{1-\frac{1}{p}}
 \left(\int_0^R\int_{{\bf R}^N}|u(x,t)|^p\psi_R^*(x,t)\,dx\,dt\right)^{\frac{1}{p}}
\end{split}
\end{equation}
for $0<R\le 1$. 
On the other hand, it follows that 
$$
\int_0^R\int_{{\bf R}^N}\chi_{\{\psi_R^*(x,t)>0\}}\,dx\,dt
=R^{\frac{N}{2m}+1}\int_0^1\int_{{\bf R}^N}\chi_{\{\psi_1^*(x,t)>0\}}\,dx\,dt.
$$
This together with \eqref{eq:3.4} implies that 
\begin{equation}
\label{eq:3.5}
\begin{split}
 & m_R+\int_0^R\int_{{\bf R}^N}|u(x,t)|^p\psi_R(x,t)\,dx\,dt\\
 & \le CR^{\frac{1}{p}\left(\frac{N(p-1)}{2m}-1\right)}
 \left(\int_0^R\int_{{\bf R}^N}|u(x,t)|^p\psi_R^*(x,t)\,dx\,dt\right)^{\frac{1}{p}}
\end{split}
\end{equation}
for $0<R\le 1$, where 
$$
m_R:=\mu\biggr(B\left(x_0,(R/3)^{\frac{1}{2m}}\right)\biggr).
$$
Let $\epsilon$ be a sufficiently small positive constant. 
For any $0<r\le R\le 1$, set 
\begin{equation}
\label{eq:3.6}
z(r):=\int_0^R\int_{{\bf R}^N}|u(x,t)|^p\psi_r^*(x,t)\,dx\,dt,
\quad
Z(R):=\int_0^R z(r)\min\{r^{-1},\epsilon^{-1}\}\,dr.
\end{equation}
Since $\eta^*$ is deceasing on $[1,\infty)$ and $\mbox{{\rm supp}}\,\eta^*\subset [1,2]$, 
for any $(x,t)\in{\bf R}^N\times(0,1)$ with $3(|x-x_0|^{2m}+t)\ge R$, 
we have 
\begin{equation}
\label{eq:3.7}
\begin{split}
 & \int_0^R \psi_r^*(x,t)\min\{r^{-1},\epsilon^{-1}\}\,dr
\le\int_0^R\eta^*\left(3\frac{|x-x_0|^{2m}+t}{r}\right)r^{-1}\,dr\\
 & \qquad\quad
\le\int_{3(|x-x_0|^{2m}+t)/R}^\infty \eta^*(s)s^{-1}\,ds\\
 & \qquad\quad
\le \eta^*\left(3\frac{|x-x_0|^{2m}+t}{R}\right)
\int_1^2 s^{-1}\,ds
\le C\psi_R^*(x,t). 
\end{split}
\end{equation}
Since $\psi_R^*(x,t)=0$ if $3(|x-x_0|^{2m}+t)<R$, 
by \eqref{eq:3.6} and \eqref{eq:3.7} we obtain 
\begin{equation}
\label{eq:3.8}
\begin{split}
 & \int_0^R\int_{{\bf R}^N}|u(x,t)|^p\psi_R(x,t)\,dx\,dt
 \ge\int_0^R\int_{{\bf R}^N}|u(x,t)|^p\psi_R^*(x,t)\,dx\,dt\\
 & \qquad\quad
\ge C^{-1}\int_0^R\int_{{\bf R}^N}|u(x,t)|^p
\left(\int_0^R \psi_r^*(x,t)\min\{r^{-1},\epsilon^{-1}\}\,dr\right)\,dx\,dt\\
 & \qquad\quad
 =C^{-1}\int_0^R\int_0^R\int_{{\bf R}^N}|u(x,t)|^p\psi_r^*(x,t)\min\{r^{-1},\epsilon^{-1}\}\,dx\,dt\,dr=C^{-1}Z(R).
\end{split}
\end{equation}
Therefore we deduce from \eqref{eq:3.5}, \eqref{eq:3.6} and \eqref{eq:3.8} that
\begin{equation}
\label{eq:3.9}
m_R+C^{-1}Z(R)
\le CR^{\frac{1}{p}\left(\frac{N(p-1)}{2m}-1\right)}(\max\{R,\epsilon\}Z'(R))^{\frac{1}{p}}.
\end{equation}
Since $m_R\ge m_r\ge m_{r_*}>0$ for $r\ge r_*$, 
it follows from \eqref{eq:3.9} that 
$$
[m_{r_*}+Z(R)]^{-p}Z'(R)\ge C^{-1}R^{-\left(\frac{N(p-1)}{2m}-1\right)}(\max\{R,\epsilon\})^{-1}
$$
for $0<r_*\le R\le 1$. 
Therefore we have 
\begin{equation}
\label{eq:3.10}
\int_{Z(r)}^{Z(1)}[m_{r_*}+s]^{-p}\,ds\ge C^{-1}\int_r^1R^{-\left(\frac{N(p-1)}{2m}-1\right)}(\max\{R,\epsilon\})^{-1}\,dR
\end{equation}
for $0<r_*\le r<1$. 
Since 
$$
\int_{Z(r)}^{Z(1)}[m_{r_*}+s]^{-p}\,ds
\le\frac{1}{p-1}(Z(r)+m_{r_*})^{-p+1}
\le\frac{1}{p-1}m_{r_*}^{-p+1},
$$
by \eqref{eq:3.10} we obtain 
$$
\frac{1}{p-1}m_{r_*}^{-p+1}
\ge C^{-1}\int_r^1R^{-\left(\frac{N(p-1)}{2m}-1\right)}(\max\{R,\epsilon\})^{-1}\,dR
$$
for $0<r_*\le r\le 1$. 
Letting $\epsilon\to+0$, 
we see that 
$$
\frac{1}{p-1}m_{r_*}^{-p+1}
\ge C^{-1}\int_r^1R^{-\frac{N(p-1)}{2m}}\,dR
$$
for $0<r_*\le r< 1$. 
This implies that
\begin{equation}
\label{eq:3.11}
\begin{split}
\mu\biggr(B\left(x_0,(r_*/3)^{\frac{1}{2m}}\right)\biggr)=m_{r_*}
 & \le C\left(\int_r^1R^{-\frac{N(p-1)}{2m}}\,dR\right)^{-\frac{1}{p-1}}\\
 & \le C\left(\int_r^{3r}R^{-\frac{N(p-1)}{2m}}\,dR\right)^{-\frac{1}{p-1}}
\le Cr^{\frac{N}{2m}-\frac{1}{p-1}}
\end{split}
\end{equation}
for $0<r_*\le r<3r<1$. 
Set $\sigma=(r/3)^{2m}=(r_*/3)^{2m}\in(0,9^{-2m})$. 
Since $x_0\in{\bf R}^N$ is arbitrary, 
we deduce from \eqref{eq:3.11} that 
\begin{equation}
\label{eq:3.12}
\sup_{x\in{\bf R}^N}\mu(B(x,\sigma))\le C\sigma^{N-\frac{2m}{p-1}},
\quad 0<\sigma<9^{-2m}.
\end{equation}
On the other hand, for any $k\ge 1$, we find $C_k>0$ such that  
\begin{equation}
\label{eq:3.13}
\sup_{x\in{\bf R}^N}\mu(B(x,k\eta))\le C_k\sup_{x\in{\bf R}^N}\mu(B(x,\eta))
\end{equation}
for $\eta>0$ (see e.g. \cite[Lemma~2.1]{IS}). 
This together with \eqref{eq:3.12} implies \eqref{eq:1.9}. 

It remains to prove \eqref{eq:1.10}. Let $p=p_m$. 
By \eqref{eq:3.11} we have 
$$
\mu\biggr(B\left(x_0,(r_*/3)^{\frac{1}{2m}}\right)\biggr)
\le C\left(\int_r^1R^{-\frac{N(p-1)}{2m}}\,dR\right)^{-\frac{1}{p-1}}
\le C|\log r|^{-\frac{N}{2m}}
\le C\biggr|\log\frac{r}{3}\biggr|^{-\frac{N}{2m}}
$$
for $0<r_*\le r<3r<1$. 
Then, similarly to \eqref{eq:3.12}, we have 
$$
\sup_{x\in{\bf R}^N}\mu(B(x,\sigma))\le C\biggr|\log\frac{r}{3}\biggr|^{-\frac{N}{2m}}
\le C\biggr[\log\left(e+\frac{1}{\sigma}\right)\biggr]^{-\frac{N}{2m}},
\quad 0<\sigma<9^{-2m}.
$$
This together with \eqref{eq:3.13} implies \eqref{eq:1.10}. 
Thus Theorem~\ref{Theorem:1.2} follows.
$\Box$
%%%%%%%%%%%%%%%%%%%%%%%%%%%%%%%%%%%%%
%%%%%%%%%%%%%%%%%%%%%%%%%%%%%%%%%%%%%
\section{Majorizing kernel}
%%%%%%%%%%%%%%%%%%%%%%%%%%%%%%%%%%%%%
%%%%%%%%%%%%%%%%%%%%%%%%%%%%%%%%%%%%%
Let $G_m=G_m(x,t)$ $(m=2,3,\dots)$ and $G_\theta=G_\theta(x,t)$ $(0<\theta<2)$
be the fundamental solutions to $\partial_t+(-\Delta)^m$ and $\partial_t+(-\Delta)^{\frac{\theta}{2}}$ 
in ${\bf R}^N\times(0,\infty)$, respectively. 
Define 
\begin{equation}
\label{eq:4.1}
K(x,t):=G_\theta\left(x,t^{\frac{\theta}{2m}}\right),\quad x\in{\bf R}^N,\,\,t>0. 
\end{equation}
Similarly to \eqref{eq:2.6} and \eqref{eq:2.7}, 
we define an integral operator $S_K(t)$ by 
$$
[S_K(t)\mu](x):=\int_{{\bf R}^N}K(x-y,t)\,d\mu(y),
\quad
[S_K(t)\phi](x):=\int_{{\bf R}^N}K(x-y,t)\phi(y)\,dy,
$$
for (signed) Radon measure $\mu$ and measurable function $\phi$ in ${\bf R}^N$. 
The aim of this section is to prove the following theorem, 
which is one of the main ingredients of this paper. 
\begin{theorem}
\label{Theorem:4.1}
Let $N\ge 1$, $m=2,3\dots$ and $\theta\in(0,2)$.  
Let $K$ be as in \eqref{eq:4.1}. 
Then $K=K(x,t)>0$ in ${\bf R}^N\times(0,\infty)$ and the following properties hold. 
\begin{itemize}
  \item[{\rm (a)}] 
  For any $j=0,1,2,\dots$, 
  there exists $d_j>0$ and $d_j'>0$ such that 
  \begin{equation*}
  |\partial_x^\alpha G_m(x,t)|\le d_jt^{-\frac{j}{2m}}K(x,t)\le d_j't^{-\frac{N}{2m}-\frac{j}{2m}}
  \end{equation*}
  for $x\in{\bf R}^N$, $t>0$ and $\alpha\in{\bf M}$ with $|\alpha|=j$. 
  \item[{\rm (b)}]
 There exists $d''>0$ such that 
 $$
 \|S_K(t)\mu\|_\infty\le d''t^{-\frac{N}{2m}}\sup_{x\in{\bf R}^N}\mu(B(x,t^{\frac{1}{2m}})),
 \quad t>0,
 $$
 for nonnegative Radon measure $\mu$ in ${\bf R}^N$. 
 \item[{\rm (c)}] There exists $d_*>0$ such that 
 \begin{equation*}
 \int_{{\bf R}^N}K(x-y,t-s)K(y,s)\,dy\le d_*K(x,t)
 \end{equation*}
 for $x\in{\bf R}^N$ and $t>s>0$.
\end{itemize}
\end{theorem}
{\bf Proof.}
The positivity of $K$ follows from the positivity of $G_\theta$ (see Section~2.2).
Let $j=0,1,2,\dots$ and $\alpha\in{\bf M}$ with $|\alpha|=j$.
By \eqref{eq:2.5} we find $C_1>0$ such that 
\begin{equation}
\label{eq:4.2}
|\partial_x^\alpha G_m(x, t)|\le C_1t^{-\frac{N}{2m}-\frac{j}{2m}}
\exp\left(-C_1^{-1}\eta^{\frac{2m}{2m-1}}\right)
\quad\mbox{with}\quad \eta=t^{-\frac{1}{2m}}|x|
\end{equation}
for $(x,t)\in{\bf R}^N\times(0,\infty)$. 
On the other hand, 
it follows from \eqref{eq:2.11} and \eqref{eq:2.12} that
$$
C_2^{-1}(1+|x|)^{-N-\theta}\le G_\theta(x,1)\le C_2(1+|x|)^{-N-\theta}, 
\quad x\in{\bf R}^N, 
$$
for some $C_2>0$. 
Then we find $C_3>0$ such that 
\begin{equation}
\label{eq:4.3}
\exp\left(-C_1^{-1}|x|^{\frac{2m}{2m-1}}\right)\le C_3G_\theta(x,1),
\quad x\in{\bf R}^N.
\end{equation}
Let $\tau:=t^{\theta/2m}$.
By \eqref{eq:2.10}, \eqref{eq:4.2} and \eqref{eq:4.3} we obtain
\begin{equation*}
\begin{split}
|\partial_x^\alpha G_m(x, t)| & 
\le C_1C_3t^{-\frac{N}{2m}-\frac{j}{2m}}G_\theta\left(t^{-\frac{1}{2m}}x,1\right)\\
 & =Ct^{-\frac{N}{2m}-\frac{j}{2m}}G_\theta\left(\tau^{-\frac{1}{\theta}}x,1\right)\\
 & =Ct^{-\frac{N}{2m}-\frac{j}{2m}}\tau^{\frac{N}{\theta}}G_\theta\left(x,\tau\right)
=Ct^{-\frac{j}{2m}}G_\theta\left(x,t^{\frac{\theta}{2m}}\right)\\
 & =Ct^{-\frac{j}{2m}}K(x,t)
\end{split}
\end{equation*}
for $(x,t)\in{\bf R}^N\times(0,\infty)$. This implies assertion~(a). 
On the other hand, 
by Lemma~\ref{Lemma:2.1} and \eqref{eq:4.1} we have 
\begin{equation*}
\begin{split}
\|S_K(t)\mu\|_\infty & =\left\|S_\theta\big(t^{\frac{\theta}{2m}}\big)\mu\right\|_\infty
\le C\big(t^{\frac{\theta}{2m}}\big)^{-\frac{N}{\theta}}
\sup_{x\in{\bf R}^N}\mu\left(B(x,\big(t^{\frac{\theta}{2m}}\big)^{\frac{1}{\theta}}\right)\\
 & =Ct^{-\frac{N}{2m}}
\sup_{x\in{\bf R}^N}\mu(B(x,t^{\frac{1}{2m}})),\quad t>0,
\end{split}
\end{equation*}
for nonnegative Radon measure $\mu$ in ${\bf R}^N$. 
This implies assertion~(b).

We prove assertion~(c). For any $0<s<t$, set 
$$
\omega_{t,s}:=(t-s)^{\frac{\theta}{2m}}+s^{\frac{\theta}{2m}}. 
$$
It follows from $\theta/2m\in(0,1)$ that 
\begin{equation}
\label{eq:4.4}
t^{\frac{\theta}{2m}}\le \omega_{t,s}=(t-s)^{\frac{\theta}{2m}}+s^{\frac{\theta}{2m}}\le 2t^{\frac{\theta}{2m}}.
\end{equation}
Then, by \eqref{eq:2.13} we have 
\begin{equation}
\label{eq:4.5}
\begin{split}
 & \int_{{\bf R}^N}K(x-y,t-s)K(y,s)\,dy
=\int_{{\bf R}^N}G_\theta(x-y,(t-s)^{\frac{\theta}{2m}})G_\theta(y,s^{\frac{\theta}{2m}})\,dy\\
 & \qquad\quad
=G_\theta(x,\omega_{t,s})
=\omega_{t,s}^{-\frac{N}{\theta}}G_\theta\left(\omega_{t,s}^{-\frac{1}{\theta}}x,1\right)
\le t^{-\frac{N}{2m}}G_\theta\left(\omega_{t,s}^{-\frac{1}{\theta}}x,1\right)
\end{split}
\end{equation}
for $x\in{\bf R}^N$ and $0<s<t$. 
Furthermore, we observe from \eqref{eq:2.11}, \eqref{eq:2.12} and \eqref{eq:4.4} that 
\begin{equation}
\label{eq:4.6}
\begin{split}
G_\theta\left(\omega_{t,s}^{-\frac{1}{\theta}}x,1\right)
 & \le C\left(1+\omega_{t,s}^{-\frac{1}{\theta}}|x|\right)^{-N-\theta}
\le C\left(1+2^{-\frac{1}{\theta}}t^{-\frac{1}{2m}}|x|\right)^{-N-\theta}\\
 & \le C\left(1+t^{-\frac{1}{2m}}|x|\right)^{-N-\theta}
\le CG_\theta\left(t^{-\frac{1}{2m}}x,1\right).
\end{split}
\end{equation}
Combining \eqref{eq:4.5} and \eqref{eq:4.6}, 
we obtain 
$$
\int_{{\bf R}^N}K(x-y,t-s)K(y,s)\,dy\le Ct^{-\frac{N}{2m}}G_\theta\left(t^{-\frac{1}{2m}}x,1\right)
=CG_\theta\left(x,t^{\frac{\theta}{2m}}\right)
=CK(x,t)
$$
for $x\in{\bf R}^N$ and $0<s<t$. 
This implies assertion~(c). Thus Theorem~\ref{Theorem:4.1} follows.
$\Box$
%%%%%%%%%%%%%%%%%%%%%%%%%%%%%%%%%%%%%
%%%%%%%%%%%%%%%%%%%%%%%%%%%%%%%%%%%%%
\section{Sufficient conditions on the solvability}
%%%%%%%%%%%%%%%%%%%%%%%%%%%%%%%%%%%%%
%%%%%%%%%%%%%%%%%%%%%%%%%%%%%%%%%%%%%
In this section, 
by use of the majorizing kernel $K$ 
we establish the existence of solutions of problem~\eqref{eq:1.1}. 
%%%%%%%%
\subsection{Existence of solutions of integral equation~(I)}
%%%%%%%%
We modify the argument in \cite{TW} to obtain 
sufficient conditions on the existence of solutions of integral equation~(I) (see Section~2.3). 
Let $T>0$ and 
$$
X:=\left\{f\in C({\bf R}^N\times(0,T))\,:\,\sup_{\tau\le t<T}\|f(t)\|_\infty<\infty\quad\mbox{for $\tau\in(0,T)$}\right\}.
$$
Let $K$ be as in Theorem~\ref{Theorem:4.1}. 
Let $U\in X$ be such that  
\begin{equation}
\label{eq:5.1}
d_*U(x,t)\ge\int_{{\bf R}^N}K(x-y,t-s)U(y,s)\,dy>0,
\quad
x\in{\bf R}^N,\,\,0<s<t<T,
\end{equation}
where $d_*$ is as in Theorem~\ref{Theorem:4.1}.
Let $\Psi$ be a positive continuous function in $(0,\infty)$ and set $V=\Psi(U)$. 
Assume that 
\begin{equation}
\label{eq:5.2}
D_*:=\sup_{0<t<T}\left\|\frac{U(t)}{\Psi(U(t))}\right\|_\infty
\int_0^t \left\|\frac{\Psi(U(s))^p}{U(s)}\right\|_\infty\,ds<\infty.
\end{equation}
Define 
$$
X_V:=\{f\in X\,:\,|||f|||<\infty\}
\quad\mbox{with}\quad
|||f|||:=\sup_{0<t<T}\sup_{x\in{\bf R}^N}\,\frac{|f(x,t)|}{V(x,t)}.
$$
Then the set $X_V$ is a Banach norm equipped with the norm $|||\cdot|||$. 
We apply the fixed point theorem in $X_V$ to prove the existence of solutions 
of integral equation~(I). 
\begin{theorem}
\label{Theorem:5.1}
Let $T>0$, $m=2,3,\dots$, $p>1$. 
Assume \eqref{eq:5.1} and \eqref{eq:5.2}.  
Let $\delta>0$ and $M>0$ be such that  
\begin{equation}
\label{eq:5.3}
\delta+d_0d_*D_*M^p\le M,\qquad 2pd_0d_*D_*M^{p-1}<1,
\end{equation}
where $d_0$ and $d_*$ are as in Theorem~{\rm\ref{Theorem:4.1}}. 
Assume that $u_0(t):=S_m(t)\mu\in X$ and $|||u_0|||\le\delta$.
Then there exists a unique solution $u\in X_V$ with $|||u|||\le M$ of integral equation~{\rm (I)} 
in ${\bf R}^N\times[0,T)$. 
\end{theorem}
{\bf Proof.}
Set 
$$
B_M:=\{u\in X_V\,:\,|||u|||\le M\}. 
$$
Define 
$$
{\mathcal F}u(t):=u_0(t)+{\mathcal N}(t),
\qquad
{\mathcal N}(t):=\int_0^t S_m(t-s)|u(s)|^p\,ds,
$$
for $u\in B_M$. Then 
\begin{equation}
\label{eq:5.4}
|{\mathcal F}u(t)|\le \delta V(t)+|{\mathcal N}(t)|,
\quad
|{\mathcal N}(t)|\le
d_0M^p\int_0^t S_K(t-s)V(s)^p\,ds.
\end{equation}
Since
\begin{equation*}
\begin{split}
 & V(x,t)^p=\frac{\Psi(U(x,t))^p}{U(x,t)}U(x,t)
\le\left\|\frac{\Psi(U(t))^p}{U(t)}\right\|_\infty U(x,t),\\
 & U(x,t)\le\frac{U(x,t)}{\Psi(U(x,t))}\Psi(U(x,t))
\le\left\|\frac{U(t)}{\Psi(U(t)}\right\|_\infty V(x,t),
\end{split}
\end{equation*}
for $(x,t)\in{\bf R}^N\times(0,T)$,  
by \eqref{eq:5.1} we have 
\begin{equation}
\label{eq:5.5}
\begin{split}
\int_0^t S_K(t-s)V(s)^p\,ds & 
\le \int_0^t \left\|\frac{\Psi(U(s))^p}{U(s)}\right\|_\infty S_K(t-s)U(s)\,ds\\
 & \le  d_*U(t)\int_0^t \left\|\frac{\Psi(U(s))^p}{U(s)}\right\|_\infty\,ds\\
 & \le d_*\left\|\frac{U(t)}{\Psi(U(t))}\right\|_\infty V(t)\int_0^t \left\|\frac{\Psi(U(s))^p}{U(s)}\right\|_\infty\,ds
\le d_*D_*V(t)
\end{split}
\end{equation}
for $0<t<T$. 
It follows from \eqref{eq:5.3}, \eqref{eq:5.4} and \eqref{eq:5.5} that 
\begin{equation}
\label{eq:5.6}
|||{\mathcal F}u|||\le \delta+d_0d_*D_*M^p\le M\quad\mbox{for}\quad u\in B_M. 
\end{equation}
On the other hand, 
by \eqref{eq:5.3} and \eqref{eq:5.5} 
we find $\nu\in(0,1)$ such that 
\begin{equation*}
\begin{split}
 & |{\mathcal F}u_1(t)-{\mathcal F}u_2(t)|\\
 & \le d_0\int_0^t S_K(t-s)||u_1|^p-|u_2|^p|\,ds\\
 & \le pd_0\int_0^tS_K(t-s)\max\{|u_1(s)|^{p-1},|u_2(s)|^{p-1}\}
V(s)\frac{|u_1(s)-u_2(s)|}{V(s)}\,ds\\
 & \le 2pd_0M^{p-1}|||u_1-u_2|||\int_0^t S_K(t-s)V(s)^p\,ds\\
 & \le 2pd_0d_*D_*M^{p-1}V(t)\,|||u_1-u_2|||
 \le\nu V(t)\,|||u_1-u_2|||
\end{split}
\end{equation*}
for $u_1$, $u_2\in B_M$. This implies that 
\begin{equation}
\label{eq:5.7}
|||{\mathcal F}u_1-{\mathcal F}u_2|||\le \nu|||u_1-u_2|||\quad\mbox{for}\quad u_1,\,u_2\in B_M. 
\end{equation}
By \eqref{eq:5.6} and \eqref{eq:5.7}
we apply the Banach fixed point theorem to find $u_*\in B_M$ uniquely such that ${\mathcal F}u_*=u_*$ in $X_V$. 
This implies that $u_*\in C({\bf R}^N\times(0,T))$ and $u_*$ satisfies
$$
u_*(x,t)=u_0(x,t)+\int_0^t\int_{{\bf R}^N}G_m(x-y,t-s)|u_*(y,s)|^p\,dy\,ds
$$
for $(x,t)\in{\bf R}^N\times(0,T)$. 
Furthermore, by \eqref{eq:5.4} and \eqref{eq:5.5} 
we have 
\begin{equation*}
\begin{split}
 & \sup_{\tau\le t<T}\|u_0(t)\|_\infty\le\delta \sup_{\tau\le t<T}V(t)<\infty,\\
 & \sup_{\tau\le t<T}\left\|\int_0^t S_m(t-s)|u_*(s)|^p\,ds\right\|_\infty
\le d_0d_*M^pD_*\sup_{\tau\le t<T}V(t)<\infty\textcolor{blue}{,}
\end{split}
\end{equation*}
for $\tau\in(0,T)$. 
Therefore we see that  $u_*$ is a solution of integral equation~(I) in ${\bf R}^N\times(0,T)$. 
Thus Theorem~\ref{Theorem:5.1} follows. 
$\Box$
%%%%%%%%%%%%%%%%%%%%
\subsection{Sufficient conditions for solvability}
%%%%%%%%%%%%%%%%%%%%
We obtain sufficient conditions for the existence of solutions of problem~\eqref{eq:1.1} 
by combining Theorem~\ref{Theorem:5.1} 
and the arguments in \cite{HI01}, \cite{RS} and \cite{W1}. 
(See also \cite{HI02}.)
We prove Theorem~\ref{Theorem:1.3}.
\vspace{5pt}
\newline
{\bf Proof of Theorem~\ref{Theorem:1.3}.}
By similar transformation~\eqref{eq:3.1} and Proposition~\ref{Proposition:2.1} 
it suffices to show the existence of solutions of integral equation~(I) in ${\bf R}^N\times[0,1)$.

We assume \eqref{eq:1.11} with $T=1$ and show the existence of solution of integral equation~(I) in ${\bf R}^N\times[0,1)$. 
Let $K$ be as in Theorem~\ref{Theorem:4.1}, that is, 
$$
K(x,t)=G_\theta\left(x,t^{\frac{\theta}{2m}}\right)\quad\mbox{with}\quad 0<\theta<2.
$$
Set $U(x,t):=2d_0[S_K(t)\mu](x)$ and $u_0(x,t):=[S_m(t)\mu](x)$. 
Then it follows from Theorem~\ref{Theorem:4.1} that 
\begin{equation*}
\begin{split}
\int_{{\bf R}^N}K(x-y,t-s)U(y,s)\,dy
 & =2d_0\int_{{\bf R}^N}\int_{{\bf R}^N}K(x-y,t-s)K(y-z,s)\,dy\,d\mu(z)\\
 & \le 2d_0d_*\int_{{\bf R}^N}K(x-z,t)\,d\mu(z)=d_*U(x,t)
\end{split}
\end{equation*}
for $x\in{\bf R}^N$ and $0<s<t$, that is, 
$U$ satisfies \eqref{eq:5.1}. 
Furthermore, it follows from Theorem~\ref{Theorem:4.1} that 
\begin{equation}
\label{eq:5.8}
|u_0(x,t)|\le d_0\int_{{\bf R}^N}K(x-y,t)\,d\mu(y)=\frac{1}{2}U(x,t),
\quad (x,t)\in{\bf R}^N\times(0,1).
\end{equation}
On the other hand,
it follows from assertion~(b) of Theorem~\ref{Theorem:4.1} and \eqref{eq:1.11} that
\begin{equation}
\label{eq:5.9}
\begin{split}
\|U(t)\|_\infty & \le Ct^{-\frac{N}{2m}}\sup_{x\in{\bf R}^N}\mu(B(x,t^{\frac{1}{2m}}))
 \le Ct^{-\frac{N}{2m}}\sup_{x\in{\bf R}^N}\mu(B(x,1))
\le C\gamma t^{-\frac{N}{2m}}
\end{split}
\end{equation}
for $0<t<1$. 
Since $1<p<1+2m/N$,  
by \eqref{eq:5.9} we have
\begin{equation}
\label{eq:5.10}
\int_0^1\|U(s)\|_\infty^{p-1}\,ds
\le (C\gamma)^{p-1}\int_0^1 s^{-\frac{N}{2m}(p-1)}\,ds
\le C\gamma^{p-1}. 
\end{equation}
We apply Theorem~\ref{Theorem:5.1} with 
$$
\Psi(s)=s,\qquad V=U,\qquad T=1,\qquad \delta=\frac{1}{2}\qquad\mbox{and}\qquad M=1.
$$
Then, by \eqref{eq:5.8} we have
\begin{equation}
\label{eq:5.11}
|||u_0|||\le\frac{1}{2}.
\end{equation}
Furthermore, by \eqref{eq:5.2} and \eqref{eq:5.10} we see that
\begin{equation}
\label{eq:5.12}
D_*\equiv\sup_{0<t\le 1}\int_0^1\|U(s)\|_\infty^{p-1}\,ds\le C\gamma^{p-1}. 
\end{equation}
Then, by \eqref{eq:5.11} and \eqref{eq:5.12}, 
taking a sufficiently small $\gamma>0$, 
we find a function $u\in B_M\subset X_V$ such that 
$$
u(t)=S_m(t)\mu+\int_0^t S_m(t-s)|u(s)|^p\,ds,\quad 0<t<1.
$$
Furthermore, we see that \eqref{eq:2.15} also holds with $T=1$. 
Therefore $u$ is a solution of integral equation~(I). 
Thus Theorem~\ref{Theorem:1.3} follows.
$\Box$
\begin{remark}
\label{Remark:5.1}
The argument in the proof of Theorem~{\rm\ref{Theorem:1.3}} is applicable 
to the case where $\mu$ is a signed Radon measure in ${\bf R}^N$. 
Indeed, the same conclusion as in Theorem~{\rm\ref{Theorem:1.3}} holds 
if $\mu$ is a signed Radon measure satisfying 
$$
\sup_{x\in{\bf R}^N}|\mu|(B(x,T^{\frac{1}{2m}}))\le\gamma_2T^{\frac{N}{2m}-\frac{1}{p-1}}
$$
for some $T>0$, instead of \eqref{eq:1.11}. 
Here $|\mu|$ is the total variation of $\mu$. 
\end{remark}
%\vspace{5pt}

Similarly to Remark~\ref{Remark:5.1}, 
we consider problem~\eqref{eq:1.1} without the nonnegativity of the initial data 
and obtain sufficient conditions for the existence of solutions of problem~\eqref{eq:1.1}. 
\begin{theorem}
\label{Theorem:5.2}
Let $N\ge 1$, $m=2,3,\dots$ and $1<\alpha<p$. 
Then there exists $\gamma=\gamma(N,m,p,\alpha)>0$ such that, 
if $\mu$ is a measurable function in ${\bf R}^N$ satisfying 
\begin{equation}
\label{eq:5.13}
\sup_{x\in{\bf R}^N}\left[\,\dashint_{B(x,\sigma)}
|\mu(y)|^\alpha\,dy\,\right]^{\frac{1}{\alpha}}\le\gamma\sigma^{-\frac{2m}{p-1}},
\qquad
0<\sigma\le T^{\frac{1}{2m}},
\end{equation}
for some $T>0$, 
then problem~\eqref{eq:1.1} possesses a solution in ${\bf R}^N\times[0,T)$. 
\end{theorem}
{\bf Proof.}
Similarly to the proof of Theorem~\ref{Theorem:1.3}, 
it suffices to show the existence of solution of integral equation~(I) in ${\bf R}^N\times[0,1)$. 
We apply Theorem~\ref{Theorem:5.1} with 
\begin{equation}
\label{eq:5.14}
\begin{split}
 & T=1,\quad u_0(x,t):=S_m(t)\mu,\quad
U(x,t):=(2d_0)^\alpha S_K(t)|\mu|^\alpha,\quad
\Psi(s):=s^{\frac{1}{\alpha}},\\
 & V(x,t):=2d_0\left(S_K(t)|\mu|^\alpha\right)^{\frac{1}{\alpha}}.
\end{split}
\end{equation}
By \eqref{eq:5.13} we have
$$
\sup_{x\in{\bf R}^N}\int_{B(x,\sigma)}|\mu(x)|^\alpha\,dx\le C\gamma^\alpha\sigma^{N-\frac{2\alpha m}{p-1}}
$$
for $0<\sigma<1$. 
This together with assertion~(b) of Theorem~\ref{Theorem:4.1} implies that 
\begin{equation}
\label{eq:5.15}
\begin{split}
\|U(t)\|_\infty & \le Ct^{-\frac{N}{2m}}
\sup_{x\in{\bf R}^N}\mu^\alpha(B(x,t^{\frac{1}{2m}}))\\
 & \le C\gamma^\alpha t^{-\frac{N}{2m}}\left(t^{\frac{1}{2m}}\right)^{N-\frac{2\alpha m}{p-1}}
\le C\gamma^\alpha t^{-\frac{\alpha}{p-1}}
\end{split}
\end{equation}
for $0<t<1$. 
By \eqref{eq:5.15} we obtain 
\begin{equation*}
\begin{split}
\int_0^t\left\|\frac{\Psi(U(s)^p}{U(s)}\right\|_\infty\,ds
 & =\int_0^t \|U(s)\|_\infty^{\frac{p-\alpha}{\alpha}}\,ds
\le C\gamma^{p-\alpha}\int_0^t s^{-\frac{p-\alpha}{p-1}}\,ds
\le C\gamma^{p-\alpha} t^{\frac{\alpha-1}{p-1}},\\
\left\|\frac{U(t)}{\Psi(U(t))}\right\|_\infty
 & =\|U(t)\|_\infty^{\frac{\alpha-1}{\alpha}}
\le C\gamma^{\alpha-1}t^{-\frac{\alpha-1}{p-1}},
\end{split}
\end{equation*}
for $0<t<1$. This implies that 
\begin{equation}
\label{eq:5.16}
D_*\equiv\sup_{0<t\le 1}\left\|\frac{U(t)}{\Psi(U(t))}\right\|_\infty
\int_0^t \left\|\frac{\Psi(U(s))^p}{U(s)}\right\|_\infty\,ds\le C\gamma^{p-1}.
\end{equation} 
On the other hand, it follows from Theorem~\ref{Theorem:4.1} and the Jensen inequality that 
\begin{equation}
\label{eq:5.17}
|S_m(t)\mu|\le d_0S_K(t)|\mu|
\le d_0\left(S_K(t)|\mu|^\alpha\right)^{\frac{1}{\alpha}}=\frac{1}{2}V(t)
\end{equation} 
for $0<t<1$. 
Similarly to the proof of Theorem~\ref{Theorem:1.3}, 
by \eqref{eq:5.16} and \eqref{eq:5.17}, 
taking a sufficiently small $\gamma>0$ and 
applying Theorem~\ref{Theorem:5.1} with \eqref{eq:5.14}, $\delta=1/2$ and $M=1$, 
we see that integral equation~(I) possesses a solution in ${\bf R}^N\times[0,1)$. 
Thus Theorem~\ref{Theorem:5.2} follows.
$\Box$
\begin{theorem}
\label{Theorem:5.3}
Let $N\ge 1$, $m=2,3,\dots$, $p=p_m$ and $\beta>0$. 
For $s>0$, set 
$$
\Phi(s):=s[\log (e+s)]^\beta,
\qquad
\rho(s):=
s^{-N}\biggr[\log\biggr(e+\frac{1}{s}\biggr)\biggr]^{-\frac{N}{2m}}. 
$$
Then there exists $\gamma=\gamma(N,m,\beta)>0$ such that,  
if $\mu$ is a nonnegative measurable function in ${\bf R}^N$ satisfying 
\begin{equation}
\label{eq:5.18}
\sup_{x\in{\bf R}^N}\Phi^{-1}\left[\,\dashint_{B(x,\sigma)}
\Phi(T^\frac{1}{p-1}|\mu(y)|)\,dy\,\right]\le\gamma\rho(\sigma T^{-\frac{1}{2m}}),
\qquad
0<\sigma\le T^{\frac{1}{\theta}},
\end{equation}
for some $T>0$, 
then problem~\eqref{eq:1.1} possesses a solution in ${\bf R}^N\times[0,T)$. 
\end{theorem}
{\bf Proof.} 
Similarly to the proof of Theorem~\ref{Theorem:1.3}, 
it suffices to show the existence of solutions of integral equation~(I) in ${\bf R}^N\times[0,1)$. 
Let $0<\gamma<1$ and assume \eqref{eq:5.18}.
Let $L\ge e$ and set $\Phi_L(s):=s[\log(L+s)]^\beta$ for $s>0$. 
Then, taking a sufficiently large $L\ge e$ if necessary, 
we have:
\begin{itemize}
  \item[{\rm (a)}] 
  $\Phi_L$ is convex in $(0,\infty)$;
  \item[{\rm (b)}] 
  the function $(0,\infty)\ni s\mapsto$
  $s^{\frac{p-1}{2}}[\log(L+s)]^{-\beta p}$ is monotone increasing.
\end{itemize}
Define a positive function $\Psi_L=\Psi_L(s)$ in $(0,\infty)$ 
by $\Psi_L(s):=\Phi_L^{-1}(s)$. 
Then 
\begin{equation}
\label{eq:5.19}
\begin{split}
C^{-1}\Phi_L(s) & \le\Phi(s)\le C\Phi_L(s),\\
C^{-1}s[\log(L+s)]^{-\beta} & \le\Psi_L(s)\le Cs[\log(L+s)]^{-\beta},
\end{split}
\end{equation}
for $s>0$. 
We apply Theorem~\ref{Theorem:5.1} with 
\begin{equation}
\label{eq:5.20}
\begin{split}
 & T=1,\quad
u_0(x,t):=S_m(t)\mu,\quad
U(x,t):=S_K(t)\Phi_L(|\mu|),\quad \Psi(s):=\Phi_L^{-1}(s),\\
 & 
V(x,t):=\Phi^{-1}_L(S_K(t)\Phi_L(|\mu|)). 
\end{split}
\end{equation}
It follows from \eqref{eq:5.18} and \eqref{eq:5.19} that 
\begin{equation}
\label{eq:5.21}
\sup_{x\in{\bf R}^N}\Phi_L^{-1}\left[\,\dashint_{B(x,\sigma)}
\Phi_L(|\mu(y)|)\,dy\,\right]\le C\gamma\rho(\sigma),
\qquad
0<\sigma<1.
\end{equation}
Applying assertion~(b) of Theorem~\ref{Theorem:4.1} with \eqref{eq:5.21}, we see that
\begin{equation}
\label{eq:5.22}
\begin{split}
\|U(t)\|_\infty & =
\|S_K(t)\Phi_L(|\mu|)\|_\infty
\le Ct^{-\frac{N}{2m}}\sup_{x\in{\bf R}^N}\int_{B(x,t^{1/2m})}
\Phi_L(|\mu(y)|)\,dy\\
 & \le Ct^{-\frac{N}{2m}}
\left(t^{\frac{1}{2m}}\right)^N
\Phi_L(C\gamma\rho(t^{\frac{1}{2m}}))\\
 & \le C\gamma\rho(t^{\frac{1}{2m}})[\log(L+C\gamma\rho(t^{\frac{1}{2m}}))]^\beta\\
 & \le C\gamma\rho(t^{\frac{1}{2m}})[\log(L+C\rho(t^{\frac{1}{2m}}))]^\beta
 \le C\gamma t^{-\frac{N}{2m}}\left|\log\frac{t}{2}\right|^{-\frac{N}{2m}+\beta}=:\gamma\xi(t)
\end{split}
\end{equation}
for $0<t<1$. 
%Let $L'\ge e$ be such that $s^{(p-1)/2}[\log(L'+s)]^{-\beta p}$ is monotone increasing in $(0,\infty)$. 
Since $p=p_m=1+2m/N$, 
it follows from property~(b), \eqref{eq:5.19} and \eqref{eq:5.22} that 
\begin{equation}
\label{eq:5.23}
\begin{split}
0 & \le\frac{\Psi_L(U(x,t))^p}{U(x,t)}
\le CU(x,t)^{p-1}[\log(L+U(x,t))]^{-\beta p}\\
 & =CU(x,t)^{\frac{p-1}{2}}U(x,t)^{\frac{p-1}{2}}[\log(L+U(x,t))]^{-\beta p}\\
 & \le C(\gamma\xi(t))^{\frac{p-1}{2}}(\gamma\xi(t))^{\frac{p-1}{2}}[\log(L+\gamma\xi(t))]^{-\beta p}\\
  & \le C\gamma^{\frac{p-1}{2}}\xi(t)^{p-1}[\log(L+\xi(t))]^{-\beta p}\\
 & \le C\gamma^{\frac{p-1}{2}}t^{-\frac{N}{2m}(p-1)}
 \left|\log\frac{t}{2}\right|^{-\frac{N}{2m}(p-1)+\beta(p-1)-\beta p}
 \le C\gamma^{\frac{p-1}{2}}t^{-1}\left|\log\frac{t}{2}\right|^{-1-\beta}
\end{split}
\end{equation}
for $(x,t)\in{\bf R}^N\times(0,1)$. 
Similarly, we have
\begin{equation}
\label{eq:5.24}
\begin{split}
0 & \le \frac{U(x,t)}{\Psi_L(U(x,t))}\le C[\log(L+U(x,t))]^\beta
\le C[\log(L+\gamma\xi(t))]^\beta\\
 & \le C[\log(L+\xi(t))]^\beta
\le C\left|\log\frac{t}{2}\right|^{\beta}
\end{split}
\end{equation}
for $(x,t)\in{\bf R}^N\times(0,1)$. 
By \eqref{eq:5.23} and \eqref{eq:5.24} we obtain 
\begin{equation}
\label{eq:5.25}
\begin{split}
D_* & \equiv\sup_{0<t<1}\left\|\frac{U(t)}{\Psi_L(U(t))}\right\|_\infty
\int_0^t \left\|\frac{\Psi_L(U(s))^p}{U(s)}\right\|_\infty\,ds\\
 & \le C\gamma^{\frac{p-1}{2}}
 \sup_{0<t<1}\biggr\{\left|\log\frac{t}{2}\right|^{\beta}\int_0^t s^{-1}\left|\log\frac{s}{2}\right|^{-1-\beta}\,ds\biggr\}
\le C\gamma^{\frac{p-1}{2}}.
\end{split}
\end{equation} 
On the other hand, it follows from Theorem~\ref{Theorem:4.1} and the Jensen inequality that 
\begin{equation}
\label{eq:5.26}
|u_0(t)|=|S(t)\mu|\le d_0S_K(t)|\mu|
\le d_0\Phi_L^{-1}\left(S_K(t)\Phi_L(|\mu|)\right)=d_0V(t)
\end{equation} 
for $0<t<1$. 
Similarly to the proof of Theorem~\ref{Theorem:1.3}, 
by \eqref{eq:5.25} and \eqref{eq:5.26}, 
taking a sufficiently small $\gamma>0$ and 
applying Theorem~\ref{Theorem:5.1} with \eqref{eq:5.20}, $\delta=d_0$ and $M=2d_0$, 
we see that integral equation~(I) possesses a solution in ${\bf R}^N\times[0,1)$.  
Thus Theorem~\ref{Theorem:5.3} follows.
$\Box$\vspace{5pt}
\newline
Theorem~\ref{Theorem:1.4} easily follows from Theorems~\ref{Theorem:5.2} and \ref{Theorem:5.3}. 
(See also \cite{HI01} and \cite{HI02}.)
\begin{remark}
\label{Remark:5.2}
Consider the Cauchy problem 
\begin{equation}
\tag{P}
\left\{
\begin{array}{ll}
\partial_t u+(-\Delta)^m u=F(u),\quad & x\in{\bf R}^N,\,\,t>0,\vspace{3pt}\\
u(x,0)=\mu(x)\ge 0, & x\in{\bf R}^N,
\end{array}
\right.
\end{equation}
where $m=2,3,\dots$ and $F$ is a continuous function in ${\bf R}$. 
Assume that 
$$
|F(u)|\le |u|^p, 
\qquad
|F(u)-F(v)|\le C_F(|u|^{p-1}+|v|^{p-1})|u-v|
$$
for $u$, $v\in{\bf R}$. 
Then, applying the same arguments in Section~{\rm 5}, 
we can show that 
the same conclusions as in Theorems~{\rm\ref{Theorem:1.3}}, {\rm\ref{Theorem:5.2}}  and {\rm\ref{Theorem:5.3}} 
and the same statement as in Remark~{\rm\ref{Remark:5.1}} hold for problem~{\rm (P)}. 
\end{remark}

\noindent
{\bf Acknowledgment.} 
The authors of this paper were supported in part 
by the Grant-in-Aid for Scientific Research (S)(No.~19H05599),
from Japan Society for the Promotion of Science.
The second author was also supported in part by the Grant-in-Aid for Young Scientists (B)
(No.~16K17629)
from Japan Society for the Promotion of Science. 
%%%%%%%%%%%%%%%%%%%%%%%%%%%%%%%%%%%%%%
%%%%%%%%%%%%    references    %%%%%%%%%%%%%%%%%%
%%%%%%%%%%%%%%%%%%%%%%%%%%%%%%%%%%%%%%

\end{document}